\documentclass{article}

\usepackage{amsthm}

\usepackage{bbm}       
\usepackage{amsmath}
\usepackage{amssymb}
\usepackage{amsfonts}
\usepackage{mathtools}

\usepackage[]{hyperref}
\usepackage{cleveref}

\usepackage{graphicx}
\usepackage{color}
\usepackage{pgfplotstable}
\usepackage{makecell}
\usepackage{booktabs}
\usepackage{colortbl}

\usepackage{enumerate}

\newcommand{\deleted}[1]{}
\newcommand{\added}[1]{{\color{black}#1}}
\newcommand{\replaced}[2]{{\color{black}#1}}

\newfont{\NUMBERS}{msbm8 scaled\magstep1}

\newcommand{\REAL}{\mathbb{R}}

\usepackage{ifthen}

\newcommand{\REALP}{\REAL_{+}}

\newcommand{\Vect}[2][]
{
  \ifthenelse{\equal{#1}{}}
  {\boldsymbol{#2}}
  {{#2}_{#1}}
}
\newcommand{\Matr}[2][]
{
  \ifthenelse{\equal{#1}{}}
  {#2}
  {{#2}_{#1}}
}

\newcommand{\Tendsto}{\rightarrow}

\newcommand{\ds}{\, \mathrm{d}s}
\newcommand{\dv}{\, \mathrm{d}V_{\Metric}}

\newcommand{\Grad}[1][]
{
  \ifthenelse{\equal{#1}{}}
  {{\nabla}}
  {{\nabla_{\!#1}}}
}
\newcommand{\Div}[1][]
{
   \ifthenelse{\equal{#1}{}}
   {{\ensuremath{\operatorname{div}}}}
   {{\ensuremath{\operatorname{div}}_{\!#1}}}
}

\def\XXint#1#2#3{{\setbox0=\hbox{$#1{#2#3}{\int}$ }
\vcenter{\hbox{$#2#3$ }}\kern-.6\wd0}}

\newcommand{\ABS}[1]{\left|#1\right|}

\newcommand{\Source}{f^{+}}
\newcommand{\Sink}{f^{-}}
\newcommand{\Forcing}{f}
\newcommand{\Supp}{\mbox{supp}}

\DeclareMathOperator{\Dist}{\operatorname{dist}}

\newcommand{\Lyap}{\mathcal{L}}
\newcommand{\Dirac}[1]{\delta_{#1}}

\newcommand{\MKEQS}{MK equations}
\newcommand{\OTP}{OT problem}
\newcommand{\OTPs}{OT problems}

\newcommand{\Cachar}{\mathcal{C}}
\newcommand{\Lip}{\operatorname{Lip}}

\newcommand{\Hilb}[1]{H^{#1}}

\newcommand{\Cont}[1][]{
	\ifthenelse{\equal{#1}{}}
	{\Cachar} 
	{\Cachar^{#1}}
}

\newcommand{\ProjMatSymb}{\mathbb{P}}
\newcommand{\ProjMat}[1][]{
  \ifthenelse{\equal{#1}{}}
  {\ProjMatSymb}
  {\ProjMatSymb_{#1}}
}
\newcommand{\IDSymbol}{\mathbb{I}}
\newcommand{\IDtens}[1][]{
  \ifthenelse{\equal{#1}{}}
  {\IDSymbol}
  {\IDSymbol(#1)}
}

\newcommand{\MPsymb}{h}
\newcommand{\MeshPar}[1][]
{
  \ifthenelse{\equal{#1}{}}
    {\MPsymb}
    {\MPsymb_{#1}}
}

\newcommand{\tstep}{k}

\newcommand{\Deltat}[1][]{\ifthenelse{\equal{#1}{}}{\Delta t_{\tstep}}{\Delta t_{#1}}}
\newcommand{\TdensH}{\Tdens_{\MeshPar}}
\newcommand{\PotH}{\Pot_{\MeshPar}}
\newcommand{\Tsymb}{\mathcal{T}}

\newcommand{\Triang}[1][]
{
  \ifthenelse{\equal{#1}{}}
  {\Tsymb}
  {\Tsymb_{#1}}
}
\newcommand{\TriangH}[1][]
{
  \ifthenelse{\equal{#1}{}}
    {\Tsymb_{\MeshPar}}
    {\Tsymb_{#1}}
}
\newcommand{\Trianghh}{\Triang_{\MeshPar/2}}

\newcommand{\Mass}{M}

\newcommand{\Test}{\phi}

\newcommand{\TolCut}{\epsilon}
\newcommand{\ACL}{\Cut{\CutPoint,\MeshPar,\TolCut}}

\newcommand{\Ene}{\mathcal{E}}
\newcommand{\Wmass}{\mathcal{M}}

\newcommand{\VspaceSymb}{\mathcal{V}}
\newcommand{\Vspace}[1][]{
  \ifthenelse{\equal{#1}{}}
  {\VspaceSymb_{\MeshPar}}
  {\VspaceSymb_{\MeshPar}(#1)}
}

\newcommand{\VbaseSymb}{\varphi}
\newcommand{\Vbase}[1][]{
  \ifthenelse{\equal{#1}{}}
  {\VbaseSymb}
  {\VbaseSymb_{#1}}
}
\newcommand{\WspaceSymb}{\mathcal{W}}
\newcommand{\Wspace}[1][]{
  \ifthenelse{\equal{#1}{}}
  {\WspaceSymb_{\MeshPar}}
  {\WspaceSymb_{\MeshPar}(#1)}
}

\newcommand{\WbaseSymb}{\psi}
\newcommand{\Wbase}[1][]{
  \ifthenelse{\equal{#1}{}}
  {\WbaseSymb}
  {\WbaseSymb_{#1}}
}

\newcommand{\SurfSymbol}{\Gamma}
\newcommand{\Surf}[1][]
{
  \ifthenelse{\equal{#1}{}}
  {{\SurfSymbol}}
  {{\SurfSymbol}_{\!#1}}
}
\newcommand{\Surfh}{\Surf[\MeshPar]}

\newcommand{\Manifold}{M}        
\newcommand{\Mdim}{n}        
\newcommand{\Chart}[1][]
{
  \ifthenelse{\equal{#1}{}}
  {\phi_{\CutPoint}}
  {\phi_{#1}}
}

\newcommand{\NormalSymbol}{\nu}
\newcommand{\Normal}[1][]{
  \ifthenelse{\equal{#1}{}}
  {{\NormalSymbol}}
  {{\NormalSymbol}_{#1}}
}

\newcommand{\Cut}[1]{C_{#1}}
\newcommand{\TangCut}[1]{\Hat{C}_{#1}}
\newcommand{\CutTime}{\mathcal{T}}
\newcommand{\TanCutTime}{\Tilde{\CutTime}}
\newcommand{\CutPoint}{p}
\newcommand{\Point}{x}
\newcommand{\SPoint}{y}

\newcommand{\Jac}{J}
\newcommand{\JacFull}{\mathcal{J}}
\newcommand{\IntSet}[1]{I_{#1}}
\newcommand{\InjDom}[1]{\Hat{I}_{#1}}
\newcommand{\Polar}{P}

\newcommand{\UnitVec}{\theta}
\newcommand{\TVect}{v}
\newcommand{\Base}{\mathcal{B}}
\newcommand{\UnitSphere}{U}
\newcommand{\Sphere}[1]{\mathcal{S}^{#1}}
\newcommand{\SphereNP}[1]{\Sphere{#1}_{\circ}}
\newcommand{\SphereDist}[1][]
{
  \ifthenelse{\equal{#1}{}}
  {{r}}
  {{r_{#1}}}
}

\newcommand{\Identification}[1][]
{
  \ifthenelse{\equal{#1}{}}
  {\Psi}
  {\Psi^{#1}}
}

\newcommand{\SphereAngle}[1][]
{
  \ifthenelse{\equal{#1}{}}
  {\varphi}
  {\varphi_{{#1}}}
}
\newcommand{\Angles}{\varphi}
\newcommand{\SphereAngles}{\mathbb{U}^{\Mdim-1}}

\newcommand{\Tdens}{\mu}
\newcommand{\Pot}{u}

\newcommand{\DMK}{DMK}
\newcommand{\OTD}{OTD}
\newcommand{\Opt}[1]{#1^*}
\newcommand{\OptTdens}{\Opt{\Tdens}}
\newcommand{\OptPot}{\Opt{\Pot}}

\newcommand{\OptTdensH}{\Opt{\TdensH}}
\newcommand{\OptPotH}{\Opt{\PotH}}

\newcommand{\Metric}[1][]
{
  \ifthenelse{\equal{#1}{}}
  {g}
  {g(#1)}
}

\newcommand{\Scal}[2]{\langle #1,#2\rangle}

\newcommand{\Norm}[2][]
{
  \ifthenelse{\equal{#1}{}}
  {{|#2|}}
  {{|#2|}_{#1}}
}
\newcommand{\NORM}[2][]
{
  \ifthenelse{\equal{#1}{}}
  {{\|#2\|}}
  {{\|#2\|}_{#1}}
}

\newcommand{\Ellipsoid}{E}
\newcommand{\Quartic}{Q}
\newcommand{\Torus}{T}

\newcommand{\Curve}{\sigma}
\newcommand{\Tangent}[2][]
{
  \ifthenelse{\equal{#1}{}}
  {T#2}
  {T_{#1}#2}
}

\newcommand{\Cellsymb}{\mathbbm{E}}
\newcommand{\Cell}[1][]{\ifthenelse{\equal{#1}{}}{\Cellsymb}{\Cellsymb_{#1}}}
\newcommand{\SubCellsymb}{\mathbbm{e}}
\newcommand{\SubCell}[1][]{\ifthenelse{\equal{#1}{}}{\SubCellsymb}{\SubCellsymb_{#1}}}
\newcommand{\CellH}[1][]{\ifthenelse{\equal{#1}{}}{\Cellsymb_{\MeshPar}}{\Cellsymb_{\MeshPar,#1}}}

\newcommand{\NCell}[1][]{
  \ifthenelse{\equal{#1}{}}
    {N_{\Cellsymb}}
    {N_{\Cellsymb({#1})}}
}

\newcommand{\Edgesymb}{\sigma}
\newcommand{\Edge}[1][]{
  \ifthenelse{\equal{#1}{}}
    {\Edgesymb}
    {\Edgesymb_{#1}}
}
\newcommand{\EdgeH}[1][]{
  \ifthenelse{\equal{#1}{}}
    {\Edgesymb_{\MeshPar}}
    {\Edgesymb_{\MeshPar,#1}}
}
\newcommand{\InradiusSymbol}{r}
\newcommand{\Inradius}[1][]
{
  \ifthenelse{\equal{#1}{}}
  {\InradiusSymbol}
  {\InradiusSymbol_{\scriptscriptstyle{#1}}}
}

\newcommand{\NeighSymbol}{\mathcal{N}}
\newcommand{\Neigh}[1][]
{
  \ifthenelse{\equal{#1}{}}
  {\NeighSymbol_{\Point}}
  {\NeighSymbol_{#1}}
}

\newcommand{\Itdcell}{r}

\newcommand{\intsphere}{\int_{\SphereAngles}}
\newcommand{\Rmax}{r_{\mbox{max}}}
\newcommand{\Rmin}{r_{\mbox{min}}}

\newcommand{\Density}{\underline{\Tdens}}
\newcommand{\IntSetTdens}{\Density}

\newcommand{\SONT}{SO(\Mdim,\Tangent[\CutPoint]{\Manifold})}
\newcommand{\SON}{SO(\Mdim)}

\renewcommand{\Tilde}[1]{#1^{\Base}}

\newtheorem{Defin}{Definition}
\newtheorem{Remark}{Remark} %

\Crefname{equation}{Eq.}{Eqs.}
\crefname{equation}{eq.}{eqs.}
\crefname{theorem}{Theorem}{Theorems}
\crefname{chapter}{Chapter}{Chapters}
\crefname{figure}{Figure}{Figures}
\crefname{Conject}{Conjecture}{Conjectures}
\crefname{Problem}{Problem}{Problems}
\crefname{Prop}{Proposition}{Propositions}        
\crefname{Theo}{Theorem}{Theorems}
\crefname{Lemma}{Lemma}{Lemma}
\crefname{Corollary}{Corollary}{Corollaries}
\crefname{section}{Section}{Sections}

\newcommand{\citet}{\cite}
\newcommand{\citep}{\cite}

\newtheorem{remark}{Remark}
\newtheorem{theorem}{Theorem}
\newtheorem{lemma}{Lemma}
\newtheorem{example}{Example}

\usepackage{PRIMEarxiv}

\usepackage[utf8]{inputenc} 

\title{Computing the Cut Locus of a Riemannian Manifold via Optimal Transport

}
\author{
  Enrico Facca\\
  Univ. Lille, Inria, CNRS, UMR 8524\\
  Laboratoire Paul Painlev\'e\\
  Lille, France \\
  \texttt{enrico.facca@inria.fr} \\
  \And
  Luca Berti\\
  CNRS\\
  Universit\'e de Strasbourg,\\
  Strasbourg, France\\
  \texttt{berti@math.unistra.fr} \\
  \And
  Francesco Fass\`o\\
  Department of Mathematics ``Tullio Levi-Civita'' \\
  University of Padova \\
  Padova, Italy\\
  \texttt{fasso@math.unipd.it} \\
  \And
  Mario Putti \\
  Department of Mathematics ``Tullio Levi-Civita'' \\
  University of Padova \\
  Padova, Italy\\
  \texttt{putti@math.unipd.it} \\
}

\begin{document}
\maketitle

\begin{abstract}
In this paper, we give a new characterization of the cut locus of a
point on a compact Riemannian manifold as the zero set of the optimal
transport density solution of the Monge-Kantorovich equations, a PDE
formulation of the optimal transport problem with cost equal to the
geodesic distance. Combining this result with an optimal transport
numerical solver based on the so-called dynamical Monge-Kantorovich
approach, we propose a novel framework for the numerical approximation
of the cut locus of a point in a manifold.  We show the applicability
of the proposed method on a few examples settled on 2d-surfaces
embedded in $\REAL^{3}$ and discuss advantages and limitations.

\end{abstract}

\newcommand{\sep}{\and}

\section{Introduction}
\label{sec:into}
Given a compact Riemannian manifold $(\Manifold,\Metric)$ of dimension $\Mdim$
and a point $\CutPoint\in \Manifold$, the cut locus $\Cut{\CutPoint}$
of $\CutPoint$ is, roughly speaking, the set of points where more than
one minimizing geodesic starting from $\CutPoint$ arrives.  For
example, the cut locus of a point in a 2-sphere embedded in $\REAL^3$
is its anti-podal point.

The cut locus is a fundamental object of Riemannian geometry, e.g., it
determines the topology of $\Manifold$ since $\Manifold\setminus
\Cut{\CutPoint}$ is diffeomorphic to an $\Mdim$-disk
~\cite{chavel_2006}.  Moreover, the cut locus is intimately related to
the singular set of the distance function~\cite{Mantegazza:2002} and
thus to the points where caustics form. However, the construction of
the cut locus of a point on a manifold is rather difficult and cut locus
shape is known only in very special cases, for example, for revolution
surfaces~\cite{Bonnard2009}.  For these reasons there is a strong
interest in its numerical construction.

The numerical approximation of the cut locus of $\CutPoint$ has seen
only sparse and diverse attempts in the past.  Many approaches are
based on the approximation of geodesics and exponential map. For
example,~\cite{Itoh-Sinclair:2004} approximate the exponential map by
means of piecewise polynomial interpolation and follow the geodesics
starting from $\CutPoint$ until they are no longer minimal. Direct
numerical evaluation of geodesics emanating from $\CutPoint$ is
another approach often proposed. This is the method of choice
in~\cite{Misztal:2011}, who use a direct discretization by finite
differences of the Riemannian partial differential equations.
In~\cite{Dey-Kuiyu:2009} the authors construct geodesics from shortest
paths of graphs constructed on point data.  All these methods require
some sort of smoothing at cut points to cope with the ill-conditioning
arising from the tracking of the distance along geodesics.  As an
alternative to geodesic-following methods, \cite{Crane:2013}
approximate the geodesic distance by means of the numerical solution
of the heat kernel defined on the manifold.  The cut locus is then
identified as the set of points where the trace of the Hessian of the
distance explodes, thus again requiring some sort of smoothing for
proper approximation.
To overcome these limitations, in~\cite{generau2020cut} the authors
propose a characterization of the cut locus as the limit in the
Hausdorff sense of a variationally-defined thawed region around the
cut locus.  This allows the construction of a convergent
finite-element-based numerical approximation of the cut locus
which is described and analyzed in~\cite{generau2020numerical}.

In this paper we propose a novel characterization of the cut locus
based on Optimal Transport (OT) theory and exploit it to derive a
stable and accurate numerical method for its approximation on compact
Riemannian manifolds, that, exploiting modern technology of numerical
methods for PDEs, is more efficient than other methods present in
literature.  The cut locus has been used extensively together with its
properties in the analysis of the regularity of optimal transport
problems~\cite{Figalli2008,Figalli:2010,Figalli-et-al:2011,Villani2011}. On
the other hand, to the best of our knowledge, our characterization of
the cut locus has never been proposed before.


In \OTPs, one looks for the optimal strategy to re-allocate a
non-negative measure $\Source$ into another non-negative measure
$\Sink$ with equal mass, given a cost for transporting one unit of
mass
(see~\cite{Ambrosio:2003,Villani:2003,Villani:2008,Santambrogio:2015}
for a complete overview of the topic).  When the transport takes place
in a Riemannian manifold with the geodesic distance as cost, the
solutions of the \OTP\ can be deduced from the solution of a nonlinear
system of PDEs known as Monge-Kantorovich equations (\MKEQS).  We
denote the solution of the \MKEQS\ by $(\OptPot,\OptTdens)$.  The
first element of the solution pair, the so-called \emph{Kantorovich
  potential} $\OptPot$, is a continuous function with Lipschitz
constant equal to $1$, whose gradient is tangent to the paths (called
rays) along which optimal transport movements occur, which are the
geodesics~\cite{Pratelli:2005}. The second solution element, the
so-called \emph{Optimal Transport Density} (\OTD) $\OptTdens$, is a
non-negative measure on $\Manifold$ that describes the mass flux
through each portion of the manifold in an optimal transportation
schedule from $\Source$ into $\Sink$.

A fundamental property of the \OTD\ $\OptTdens$ is that it decays
towards zero at the endpoints of the geodesics along which the mass is
moved~\cite{Evans-Gangbo:1999}. This fact suggests that if we take
$\Source=\dv(\Manifold)\Dirac{\CutPoint}$ and $\Sink=\dv$ (where
$\Dirac{\CutPoint}$ denotes the Dirac measure centered at $\CutPoint$,
$\dv$ the volume form induced by the metric $\Metric$, and
$\dv(\Manifold)$ is the measure of $\Manifold$), the \OTD\ restricted
to all geodesics starting at point $\CutPoint$ tends to zero when
approaching the points that form the cut locus of $\CutPoint$ in
$\Manifold$. This intuition is confirmed by the following theorem,
which represents our characterization of $\Cut{\CutPoint}$:
\begin{theorem}
  \label{thm:opts-intro}
  Let $(\Manifold,\Metric)$ be a compact and geodesically complete
  Riemannian manifold of dimension $\Mdim$ with no boundary.
  Given a point $\CutPoint\in\Manifold$, the \OTD\ $\Opt{\Tdens}$
  solution of the \MKEQS\ with
  $\Source=\dv(\Manifold)\Dirac{\CutPoint}$ and $\Sink=\dv$ admits, in
  the set $\Manifold \setminus \left\{\CutPoint\right\}$, a continuous
  density $\Density$ with respect to the volume form $\dv$ whose zero
  set coincides with $\Cut{\CutPoint}$ i.e.,
  \begin{equation*}
    \Opt{\Tdens} = \Density \dv 
    \ , \ 
    \Cut{\CutPoint}
    =
    \left\{
      \Point \in \Manifold \setminus \left\{\CutPoint\right\}
      \ : \ \Density(\Point) = 0  
    \right\}
    .
  \end{equation*}  
\end{theorem}

Based on this characterization, we propose a new
\added{variationally-based} numerical scheme to approximate the cut
locus on compact surfaces by means of the numerical solution of the
\MKEQS.  For the latter, we adopt the approach described
in~\cite{Berti-et-al-numerics:2021}, where the discrete Dynamical
Monge-Kantorovich (DMK) framework described
in~\cite{Facca-et-al-numeric:2020,Facca-et-al:2018} is extended to
$\REAL^3$-embedded surfaces.  We show the effectiveness of the
proposed numerical approach identifying the cut locus of the following
triangulated surfaces: a torus, for which the cut locus is known, and
two test cases borrowed from~\cite{Itoh-Sinclair:2004} to show the
applicability to generic surfaces.

The paper is organized as follows. First we present all the Riemannian
objects required for the definition of the cut-locus, together with
some of its properties.  Then, in~\cref{sec:otp-manifold} we recall
the definition of the Optimal Transport Problem with cost equal to the
geodesic distance and the \MKEQS. \Cref{sec:cut-locus-otp} is
dedicated to the connection of the \MKEQS\ with the cut-locus and ends
with the proof of~\cref{thm:opts-intro}. Finally, \Cref{sec:cut
  locus-experiments} is dedicated to the presentation of the proposed
numerical approach to identify the cut locus of triangulated surfaces
via the DMK strategy and concludes with some numerical experiments.

\section{The cut locus of a point on a Riemannian Manifold}
\label{sec:manifold}
We consider a geodesically complete and compact Riemannian manifold
$(\Manifold,\Metric)$ with no boundary, equipped with a smooth metric
$\Metric$.  We denote with $\Scal{v}{w}_{\Metric[\Point]}$ the
application of the metric $\Metric$ evaluated at $\Point$ to the two
vectors $v,w\in\Tangent[\Point]{\Manifold}$.  The symbols
$\Grad[\Metric]$, $\Norm[\Metric]{\cdot}$, and $\dv$ are used to
identify the gradient operator, the vector norm, and the volume form
induced by the metric tensor $\Metric$ on $\Manifold$.  The distance
between two points $\Point,y\in \Manifold$ is defined as:
\begin{equation*}
 \Dist_{\Metric}(x,y) = \inf_{\Curve}
 \left\{
   \int_{0}^{1}
   \sqrt{
     \Scal{\dot{\Curve}(s)}{\dot{\Curve}(s)}_{\Metric[\Curve(s)]}
   } 
   \; \mathrm{d}s
   \;:\; 
   \begin{gathered}
     \Curve\in \Cont[1]([0,1],\Manifold)\\
     \Curve(0) = x, \Curve(1) =y
   \end{gathered}
 \right\}.
\end{equation*}
Given a point $\CutPoint\in \Manifold$ , we denote with
$\exp_{\CutPoint}:\Tangent_{\CutPoint}{\Manifold}\rightarrow\Manifold$
the exponential map of $(\Manifold,\Metric)$ at $\CutPoint$ 
and with $\exp^{-1}_{\CutPoint}$ its inverse, where well defined.

\subsection{Cut locus of a point}
\label{sec:cut locus}

\begin{figure}
  \centerline{
    \includegraphics[width=0.8\textwidth,
    trim={0.0cm 0.6cm 1cm 0.0cm},clip
    ]{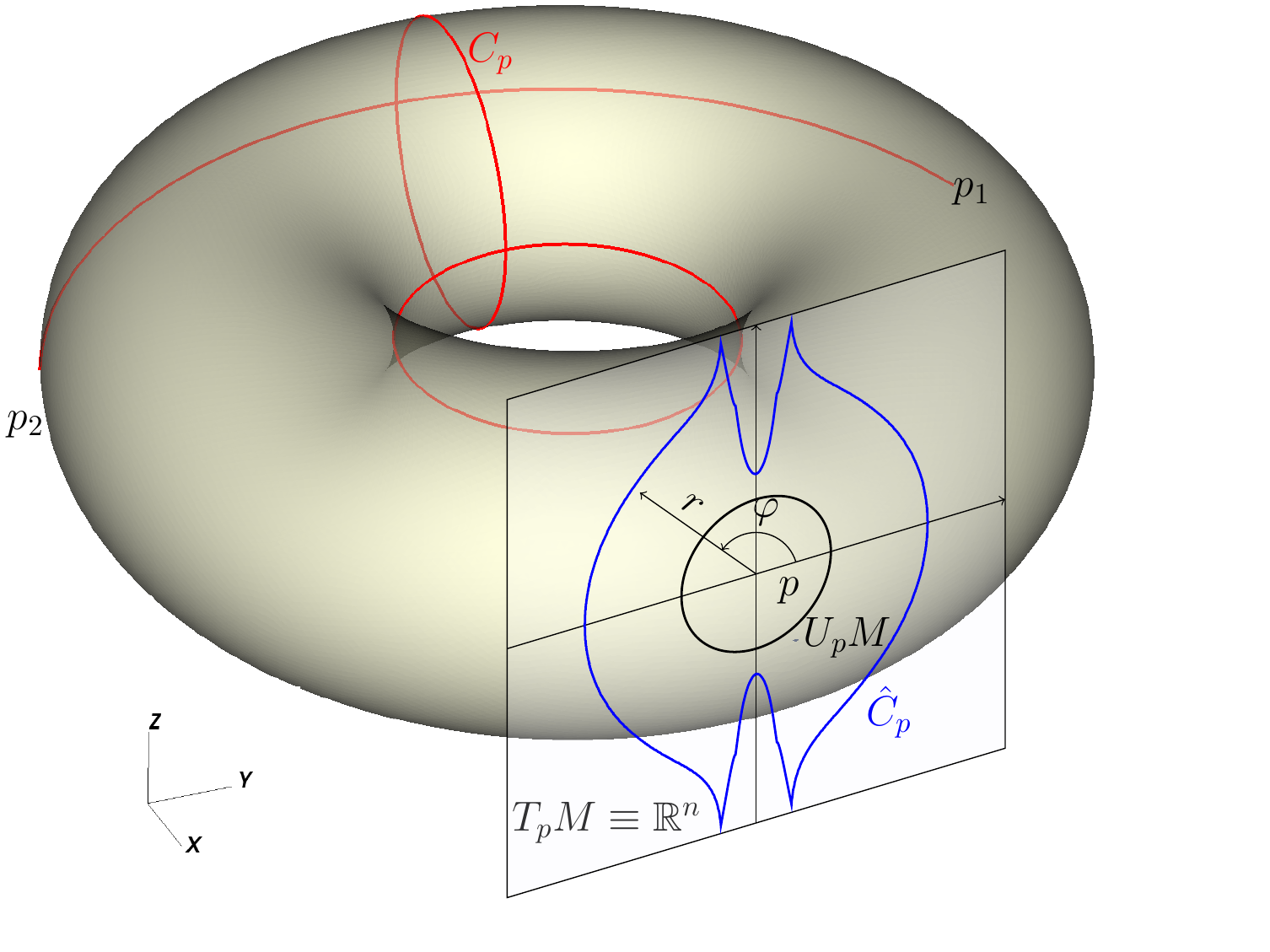}
  }
  \caption{The tangent cut locus and the cut locus of a point on the
    torus with major and minor radi $\Rmax=2$ and $\Rmin=1$ embedded
    in $\REAL^{3}$. The red curves identify the cut locus of the point
    $\CutPoint=(R,0,0)$ with $R:=\Rmax+\Rmin$ located on the external
    ``equator''. Note that the cut locus is formed by the internal
    ``equator'', the opposite ``meridian'', and a portion of the
    external ``equator'' connecting the two points
    $p_1=(R\cos(\alpha),R\sin(\alpha),0)$ and
    $p_2=(R\cos(\alpha),-R\sin(\alpha),0)$ with \replaced{$\alpha =
      \pi\Rmin /\sqrt{R \Rmin}$}{$\alpha = \pi/\sqrt{R \Rmin^2}$}
    (formulae taken from~\cite[Theorem 9]{kimball1930}). The blue
    lines represent\deleted{s} \replaced{an approximation of the
      tangent cut locus $\TangCut{\CutPoint}$ obtained by
      straight-forward numerical approximation of the geodesics
      emanating from $p$.}{the piecewise linear interpolation of the
    tangent cut locus through some meaningful points of $\Manifold$,
    namely $p_1,p_2$ and points $(\Rmin,0,0),(-\Rmin,0,0),(-R,0,0)$.}
  }
  \label{fig:torus-chart}
\end{figure}

We now give the definition and some properties of the cut locus
$\Cut{\CutPoint}$ of a point $\CutPoint\in\Manifold$ and of all the
related objects that will be used in the sequel
(see~\cite{Sakai1996}).  \Cref{fig:torus-chart} illustrates
graphically these definitions for a torus embedded in $\REAL^{3}$.

\begin{Defin}
  \label{def:manifold}
  Let $(\Manifold,\Metric)$ be a compact and geodesically complete Riemannian
  manifold of dimension $\Mdim$ and consider a point
  $\CutPoint\in\Manifold$. Let $\UnitSphere_{\CutPoint}\Manifold$ be
  the set of unit tangent vectors at $\CutPoint$, i.e.:
  \begin{equation*}
    \UnitSphere_{\CutPoint}\Manifold
    := 
    \left\{
      \UnitVec\in\Tangent[\CutPoint]{\Manifold }
      \ : \Scal{\UnitVec}{\UnitVec}_{\Metric[\CutPoint]}=1
    \right\} \; .
  \end{equation*}
  The \emph{cut time} $\CutTime_{\CutPoint}(\UnitVec)$ of
  $\UnitVec\in \UnitSphere_{\CutPoint}\Manifold$ is defined as:
  \begin{align}
    \label{eq:cut-time}
    \CutTime_{\CutPoint}(\UnitVec) := \sup\left\{ t\geq 0 \,: \,
      \exp_{\CutPoint}(s\UnitVec)_{0\leq s< t}\  \mbox{is a
        minimizing geodesic} \right\} \; .
  \end{align}
  The sets:
  \begin{equation}
    \InjDom{\CutPoint}:=
    \left\{
    t \UnitVec
    \, : \,
    \UnitVec \in \UnitSphere_{\CutPoint}\Manifold
    \, , \,
    0< t<\CutTime_{\CutPoint}(\UnitVec) 
    \right\}
    \subset
    \Tangent[\CutPoint]{\Manifold }
    \quad\mbox{and}\quad
    \IntSet{\CutPoint}:=\exp_{\CutPoint}(\InjDom{\CutPoint})
    \subset \Manifold
  \end{equation}
  are called the \emph{injectivity domain} and the \emph{interior set}
  at $\CutPoint$.  Note that $\exp_{\CutPoint}$ is a diffeomorphism from
  $\InjDom{\CutPoint}$ to $\IntSet{\CutPoint}$, thus
  $\exp_{\CutPoint}^{-1}$ is defined from $\IntSet{\CutPoint}$ to
  $\InjDom{\CutPoint}$.\\
  The sets:
  \begin{equation*}
    \TangCut{\CutPoint}
    :=
    \left\{
    \CutTime_{\CutPoint}(\UnitVec) \UnitVec 
    \, : \,
    \UnitVec \in \UnitSphere_{\CutPoint}{\Manifold}
    \right\}
    \subset  \Tangent[\CutPoint]{\Manifold}
    \quad\mbox{and}\quad
    \Cut{\CutPoint}:=\exp_{\CutPoint}(\TangCut{\CutPoint})
    \subset \Manifold
  \end{equation*}
  are called, respectively, the \emph{tangent cut locus} and the
  \emph{cut locus} of $\CutPoint$.
\end{Defin}
\begin{lemma}[Properties of the cut time \cite{Itoh-Minoru:2001}]
  \label{lemma:cut-time-prop}
  The functions $\CutTime_{\CutPoint}$ and $\TanCutTime_{\CutPoint}$
  defined in~\cref{eq:cut-time,eq:cut-time-tilde} are Lipschitz
  continuous.
\end{lemma}
\begin{lemma}[Properties of the cut locus~\cite{Sakai1996}]
  \label{lemma:cut-locus-prop} 
  The manifold $\Manifold$ is the union of three disjoint sets: the
  point $\CutPoint$, the interior set $\IntSet{\CutPoint}$, and the
  cut locus $\Cut{\CutPoint}$.  Moreover, $\Cut{\CutPoint}$ is a null
  set of the volume form $\dv$.
\end{lemma}

\begin{lemma}[Properties of the injectivity domain~\cite{Sakai1996}]
  \label{lemma:inj-prop} 
 The injectivity domain $\InjDom{\CutPoint}$ is a
 star shaped subset of $\Tangent[\CutPoint]{\Manifold}$ with boundary
 given by $\TangCut{\CutPoint}$.
\end{lemma}

As an example,~\cref{fig:torus-chart} shows the cut locus
$\Cut{\CutPoint}$ (red color) of a point $\CutPoint$ lying
on the external ``equator'' of the torus (from~\cite{Gravensen2005,kimball1930})
and  \deleted{a piecewise linear interpolation of some meaningful
  points of}
\added{an approximation of } the
tangent cut locus
$\TangCut{\CutPoint}\subset\Tangent[\CutPoint]{\Manifold}$ (blue
color) \added{obtained by solving numerically the equations of
  the geodesic curves emanating from $p=(3,0,0)$}.

\subsection{Riemannian polar coordinates}
\label{sec:polar-coordinates}
Next we give workable expressions, in local coordinates, for
these cut locus related quantities. To this aim, we need to fix proper
charts and local coordinates. In view of~\cref{lemma:inj-prop}, it is
convenient to use Riemannian polar coordinates (see~\cite{Jost:2008}\added{)}.
However we need to be slightly pedantic here about their definition
and properties because we will need to show that our solution of the OT
problem is free of the singularities introduced by such coordinate
systems).

Any set of these coordinates depends on the choice of a
$\Metric[\CutPoint]$-orthonormal basis
$\Base=\left\{e_1^{\Base},\ldots,e_\Mdim^{\Base}\right\}$ of
$\Tangent[\CutPoint]{\Manifold}$, which allows the identification of
$\Tangent[\CutPoint]{\Manifold}$ with $\REAL^{\Mdim}$ via the map
$\TVect\mapsto (\TVect^{\Base}_1,\ldots,\TVect^{\Base}_\Mdim)$ with
$\TVect^{\Base}_i:=\Scal{\TVect}{e_i^{\Base}}_{\Metric[\CutPoint]}$
for $i=1,\ldots,\Mdim$.  We define the map
\begin{align*}
  &\Identification[\Base]:\InjDom{\CutPoint}\to \REALP\times \Sphere{\Mdim-1}
  \ , \ 
  \TVect \mapsto \left(\Norm[\Metric]{\TVect}, \frac{\TVect^{\Base}}{\Norm[\Metric]{\TVect}}\right)
\end{align*}
and, denoting $\SphereNP{\Mdim-1}=\Sphere{\Mdim-1}\setminus
(0,\ldots,\pm 1)$, we introduce the (Euclidean) polar coordinates
\begin{align*}
  \Polar: \REALP \times \SphereNP{\Mdim-1}
  &\to
  \REALP \times \SphereAngles
  \quad , \quad 
  \SphereAngles
  = ]0,\pi[^{\Mdim-2}\times \Sphere{1}
      \\
      (\SphereDist,\SPoint)&\mapsto
      \left(\SphereDist, \Angles_1(\SPoint),\ldots,\Angles_{\Mdim-1}(\SPoint) \right)
\end{align*}
where $\Angles=\left(\Angles_1,\ldots, \Angles_{\Mdim-1} \right)$ are the
$\Mdim-1$ spherical coordinates parameterizing $\Sphere{\Mdim-1}$
(see~\cite{Blumenson:1960} for explicit formulas).
Then, the Riemannian polar coordinates relative to the basis $\Base$
is the map $\Tilde{\Chart}= \Polar \circ \Identification[\Base] \circ
\exp_{\CutPoint}$. Explicitly, they are given by:
\begin{alignat}{4}
  \label{eq:tilde-chart}
  \Tilde{\Chart}
  :\IntSet{\CutPoint}^{\Base}
  &\overset{\exp_{\CutPoint}}{\longrightarrow}
  \ \ \  \InjDom{\CutPoint}^{\Base}
  &&\overset{\Identification[\Base]}{\longrightarrow}
  \REALP \times \SphereNP{\Mdim-1}
  &&\overset{\Polar}{\longrightarrow}
  \REALP \times \SphereAngles
  \\
  \nonumber
  \Point &\xmapsto{\ \ }
  \TVect:=\exp_{\CutPoint}(\Point)
  &&\xmapsto{\ \ }
  \left(\Norm[\Metric]{\TVect}
  ,
  \frac{\TVect^{\Base}}{\Norm[\Metric]{\TVect}}\right)
  &&\xmapsto{\ \ }
  \left(\Norm[\Metric]{\TVect},
  \SphereAngle\left(\frac{\TVect^{\Base}}{\Norm[\Metric]{\TVect}}\right)\right)
\end{alignat}
where the sets $\IntSet{\CutPoint}^{\Base}\subset\IntSet{\CutPoint}$
and $\InjDom{\CutPoint}^{\Base}\subset \InjDom{\CutPoint}$ are the
preimages under the maps $\Identification[\Base] \circ
\exp_{\CutPoint}$ and $\Identification[\Base]$ of $\REALP \times
\SphereNP{\Mdim-1}$, respectively. Note that
$\Norm[\Metric]{\TVect}=\Dist_{\Metric}(\CutPoint,\Point)$.

\begin{Remark}
  \label{remark:line}
  The set $\IntSet{\CutPoint}^{\Base}$ coincides with
  $\IntSet{\CutPoint}$ minus a geodesic curve which passes through
  $\CutPoint$ and depends on the basis $\Base$.
\end{Remark}

It is clear that, given $\Base$, there exists a one-to-one
correspondence between elements of $\UnitSphere_{\CutPoint}\Manifold\setminus\left\{e_{\Mdim}^{\Base},-e_{\Mdim}^{\Base} \right\}$
and $\SphereAngles$.  Thus, we define the cut time function
$\TanCutTime_{\CutPoint}$ that maps the angle variables $\SphereAngle$
into the cut time of the corresponding unit vector in the tangent
space i.e.,
\begin{equation}
  \label{eq:cut-time-tilde}
  \TanCutTime_{\CutPoint}:\SphereAngles \to \REALP
  \ , \ 
  \TanCutTime_{\CutPoint}(\SphereAngle) :=
  \CutTime_{\CutPoint} \circ \left(\Identification[\Base]\right)^{-1}\circ
  \Polar^{-1}(1,\SphereAngle) \quad \forall \SphereAngle\in
  \SphereAngles
\end{equation}

\begin{lemma}[Properties of polar coordinates,~\cite{Jost:2008}]
  \label{lemma:metric-polar-coordinates}
  For any $\Metric[\CutPoint]$-orthonormal basis $\Base$, the matrix
  representing the metric $\Metric|_{\IntSet{\CutPoint}^{\Base}}$
  written in Riemannian polar coordinates has the block-diagonal
  expression
  \begin{gather*}
    \Tilde{\Metric}(\SphereDist,\SphereAngle)
    =
    \begin{pmatrix}
      1 & 0 \\
      0 & \Tilde{h}(\SphereDist,\SphereAngle)
    \end{pmatrix}
    \; ,
  \end{gather*}
  where $\Tilde{h}(\SphereDist,\SphereAngle)$ is a
  $(\Mdim-1)\times(\Mdim-1)$ symmetric and
  positive-definite matrix.
\end{lemma}
We denote
\begin{equation*}
  \JacFull^{\Base}(\SphereDist,\SphereAngle)d \SphereDist \wedge d\SphereAngle[1]\wedge\ldots\wedge d\SphereAngle[\Mdim-1]
\end{equation*}
the volume form $\dv$ expressed in polar coordinates in the interior
set $\IntSet{\CutPoint}$. 
\begin{lemma}
  \label{lemma:factorization}
  The function $\JacFull^{\Base}$ factorizes as follows
  \begin{equation}
    \label{eq:polar-volume-form}
    \JacFull^{\Base}(\SphereDist,\SphereAngle) 
    =
    G^{\Base}(\SphereDist,\SphereAngle )
    \Jac( \SphereDist,\SphereAngle ) 
    \; ,
  \end{equation}
  where $G^{\Base}:\REALP\times \SphereAngles\to \REALP$ is given by
  \begin{equation}
    \label{eq:metric-in-coordinate}
    G^{\Base}(\SphereDist,\SphereAngle) 
    =
    \sqrt{
      \det
      \left(
      \Metric \circ (\Chart^{\Base})^{-1}( \SphereDist,\SphereAngle )
      \right)
    }
  \end{equation}
  and $\Jac$ is the absolute value of the determinant of the Jacobian
  matrix of $\Polar^{-1}$ and is given by:
  \begin{gather}
  \label{eq:polar-jacobian}
  \Jac(\SphereDist,\SphereAngle^{1},\SphereAngle^{2},\ldots,\SphereAngle^{\Mdim-2}) = 
    \SphereDist^{\Mdim-1}
    \sin^{\Mdim-2}(\SphereAngle_{1})
    \sin^{\Mdim-3}(\SphereAngle_{2})\ldots
    \sin(\SphereAngle_{\Mdim-2})
     \; .
  \end{gather}
\end{lemma}
\begin{proof}
  The proof comes directly from the appropriate composition of the
  relevant functions defined above.
\end{proof}

Let $\SONT$ be the set of all proper ($\det=+1$) linear
$\Metric(\CutPoint)$-isometries of $\Tangent[\CutPoint]{\Manifold}$.
Given a $\Metric(\CutPoint)$-orthonormal basis
$\Base=\left\{e_1^{\Base},\ldots,e_\Mdim^{\Base}\right\}$ of
$\Tangent[\CutPoint]{\Manifold}$, we associate to any
$\Matr{\hat{R}}\in \SONT$ the matrix $\Matr{R}^{\Base}\in \SON$ with
entries
$\Matr[i,j]{R}^{\Base}=
\Scal{e^{\Base}}{\Matr{\hat{R}}e^{\Base}_j}_{\Metric[\CutPoint]}$. 
Also, we denote by $\Matr{\hat{R}}\Base$ the rotated basis
$\left\{ \Matr{\hat{R}} e^{\Base}_1, \ldots, \Matr{\hat{R}}
  e^{\Base}_\Mdim\right\}$. 

\begin{lemma}
  \label{lemma:charts}
  For any $\Metric(\CutPoint)$-orthonormal basis $\Base$ and any
  $\Matr{\hat{R}}\in \SONT$,
  \begin{equation*}
  G^{\Matr{\hat{R}} \Base}\circ \Polar(\SphereDist,\Matr{R}^{\Base} \SPoint)
  =
  G^{\Base}\circ \Polar(\SphereDist, \SPoint)
  \end{equation*}
  for all $\SphereDist \in \REALP$ and $\SPoint \in\SphereNP{\Mdim-1}$
  such that $(\Matr{\hat{R}})^T\SPoint\in \SphereNP{\Mdim-1}$. 
\end{lemma}

\begin{proof}
  This follows from the fact that, as we verify,
  \begin{equation}
    \label{eq:theo-stat}
    \left(
    \Identification[\Matr{\hat{R}} \Base]
    \right)^{-1}
    (\SphereDist, \Matr{R}^{\Base} \SPoint)
    =
    \left(
    \Identification[\Base]
    \right)^{-1}
    (\SphereDist,\SPoint )
  \end{equation}
  for all $\SphereDist, \SPoint$ as in the statement.  Write for short
  $\Matr{R}$ for $\Matr{R}^{\Base}$.  First note that, $\forall \TVect
  \in \Tangent[\CutPoint]{\Manifold}$,
  $\TVect^{\Matr{\hat{R}}\Base}=\Matr{R}^T\TVect^{\Base}$ and
  $(\Matr{\hat{R}}\TVect)^{\Base}=\Matr{R}\TVect^{\Base}$.
  Thus
  \begin{align*}
    \Identification[\Matr{\hat{R}} \Base](\Matr{\hat{R}}\TVect)
    &=\left(
    |\Matr{\hat{R}}\TVect|_{\Metric} ,
    \frac{(\Matr{\hat{R}}\TVect)^{\Matr{\hat{R}}\Base}}{|\Matr{\hat{R}}\TVect|_{\Metric}}
    \right)
    =\left(
    |\TVect|_{\Metric} ,
    \frac{\Matr{R}^T (\Matr{\hat{R}}\TVect)^{\Base}}{|\TVect|_{\Metric}}
    \right)\\
    &=\left(
    |\TVect|_{\Metric} ,
    \frac{\Matr{R}^T \Matr{R}\TVect^{\Base}}{|\TVect|_{\Metric}}
    \right)
    = \Identification[\Base](\TVect)\; .
  \end{align*}
  This implies that, if $(\SphereDist, \SPoint)= \Identification[\Base](\TVect)$
  then
  \begin{align*}
    \Identification[\Matr{\hat{R}} \Base] \circ
    \left(\Identification[\Base]\right)^{-1}(\SphereDist, \SPoint)
    &=
    \Identification[\Matr{\hat{R}} \Base](\TVect)
    =
    \left(
    |\TVect|_{\Metric} ,
    \frac{\TVect^{\Matr{\hat{R}} \Base}}{|\TVect|_{\Metric}}
    \right)
    \\
    &= \left(
    |\TVect|_{\Metric} ,
    \frac{\Matr{R}\TVect^{\Base}}{|\TVect|_{\Metric}}
    \right)
    =
    \left(
    \SphereDist, \Matr{R} \SPoint
    \right)
    \; ,
  \end{align*}
  namely~\cref{eq:theo-stat}.
\end{proof}

\section{Monge-Kantorovich equations on manifolds}
\label{sec:otp-manifold}
In this section we present the Monge-Kantorovich equations
(\MKEQS), an equivalent PDE formulation of the
optimal transport problem on a Riemannian manifold with geodesic
distance as transport cost.  We use the formulation described
in~\cite{Bouchitte-Buttazzo:2001,Pratelli:2005}, assuming that one of
the transported measures $\Source$ and $\Sink$ on $\Manifold$ admits a
density with respect to the volume form $\dv$.  Under these
assumptions the \MKEQS\ can be written as the problem of finding a
pair $(\OptPot,\OptTdens)$, where $\OptPot$ is a continuous function
with Lipschitz constant equal to $1$ and $\OptTdens$ is a non-negative
measure\deleted{s}, that solves
\begin{subequations}
  \label{eq:mkeqs-manifold}
  \begin{align}
    \label{eq:mkeqs-manifold-elliptic}
    -\Div[\Metric] (\OptTdens \Grad[\Metric] \OptPot)
    &= 
    \Source -\Sink 
    \quad 
    &\textrm{on $\Manifold$}
    \\
    \label{eq:mkeqs-manifold-gradless1}
    \Norm[\Metric]{\Grad[\Metric]\OptPot} &\leq 1 
    \quad &\textrm{on $\Manifold$}
    \\
    \label{eq:mkeqs-manifold-eikonal}
    \Norm[\Metric]{\Grad[\Metric]\OptPot} &= 1 \quad
    &\OptTdens-a.e.
  \end{align}
\end{subequations}
where \cref{eq:mkeqs-manifold-elliptic} must be interpreted in 
the following weak form:
\begin{equation}
  \label{eq:mkeqs-weak-manifold}
  \int_{\Manifold} 
  \Scal{
    \Grad[\Metric]\OptPot
  }{
    \Grad[\Metric]\Test}_{\Metric}
  \mathrm{d}\OptTdens
  =
  \int_{\Manifold}
  \Test
  \mathrm{d}\Source
  - 
  \int_{\Manifold}
  \Test
  \mathrm{d}\Sink
  \ , \ 
  \forall \Test \in \Cont[1](\Manifold) \; .
\end{equation}
The components of the solution pair $(\OptPot,\OptTdens)$ of the above
system are named Kantorovich potential and \OTD, respectively.

The following lemma summarizes a series of results on the solution of
the
\MKEQS\ from~\cite{Feldman-McCann:2001,Feldman-McCann:2002,Ambrosio:2003,DePascale-et-al:2004,Santambrogio:2009,Pratelli:2005}.
%
\begin{lemma}[Properties of Monge-Kantorovich equations]
  \label{lemma:prop-tdens}
  The Kantorovich potential $\OptPot$ is unique up to a constant
  within the support of $\OptTdens$ (outside $\Supp(\OptTdens)$ there
  exist infinitely many functions $\Pot$ that
  satisfy~\cref{eq:mkeqs-manifold}~\cite{Evans-Gangbo:1999}).
  
  If either $\Source$ or $\Sink$ (or both) is absolutely continuous
  with respect to the volume form $\dv$, then the \OTD\ $\OptTdens$
  is unique and is absolutely continuous
  with respect to $\dv$.
\end{lemma}
\begin{Remark}
  Without a priori knowledge of the integrability of $\OptTdens$,
  guaranteed by~\cref{lemma:prop-tdens} in our problem, the formulation
  of the \MKEQS\ requires the notion of gradient with respect to
  measure~\cite{Bouchitte-Buttazzo:2001}.
\end{Remark}

Intuitively speaking, the link between OT and the solution
$(\OptPot,\OptTdens)$ of the \MKEQS~\cref{eq:mkeqs-manifold} is as
follows. The \OTD\ $\OptTdens$ can be seen as a measure of the flux
through each portion of the manifold in the optimal reallocation of
$\Source$ into $\Sink$.  Mass moves along disjoint transport rays that
follow the direction of the gradient of the Kantorovich potential
$\OptPot$. Intuitively, these transport rays are geodesics connecting
points in the support of $\Source$ with points in the support of
$\Sink$ (we refer the reader to~\cite{Ambrosio:2003} for the proper
definition of transport rays). Under certain properties of the
transported measures $\Source$ and $\Sink$ ($L^\infty$ densities and
disjoint supports),~\cite{Evans-Gangbo:1999} proved that the
\OTD\ tends to zero at the endpoints of each transport ray.

\section{OT characterization of $\Cut{\CutPoint}$}
\label{sec:cut-locus-otp}
From the property of the decay of the \OTD\ along the transport rays
we devise the following strategy for the search of the cut locus: we
set $\Source=\dv(\Manifold)\Dirac{\CutPoint}$ (with
$\Dirac{\CutPoint}$ the Dirac delta centered at $\CutPoint$) and
$\Sink=\dv$ and look at the zero-set of $\OptTdens$.  The intuition is
that in the optimal reallocation of the Dirac mass centered at
$\CutPoint$, the mass is ``sent'' from $\CutPoint$ to all the points
of $\Manifold$ along geodesics. The mass is progressively ``absorbed''
by the constant sink term $\Sink=\dv$ until we reach those points
where mass is coming also from a different direction.  At those points,
$\OptTdens$ becomes zero. For a cleaner understanding of this idea and
of the properties of the \OTD, we present a simple example on the unit
circle $\mathcal{S}^1$.
\begin{example}
  \begin{figure}
    \centerline{
      \includegraphics[width=0.7
        \textwidth,
        trim={3cm 22.50cm 12cm 1.0cm},clip,
      ]{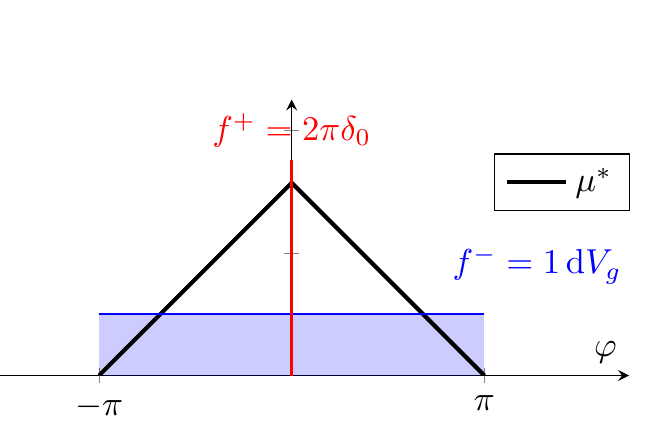}
    }
    \caption{ Graph of the \OTD\ $\OptTdens$ (black) solution of
      the \MKEQS\ in the unit circle $\mathcal{S}^1$
      (the point $-\pi$ is identified with $\pi$)
      with
      $\Source=2\pi\Dirac{0}$ (red) and
      $\Sink=d\Angles$ (blue)}
    \label{fig:circle}
  \end{figure}
  On the unit circle $\mathcal{S}^1=\left\{x^2+y^2=1\right\}$ we
  consider the point $\CutPoint=(1,0)$ and source terms given by
  $\Source=2\pi\Dirac{\CutPoint}$ and $\Sink=1d\Angles$. Using the
  angle coordinate $\Angles\in \Sphere{1}$ the Kantorovich potential
  and the \OTD\ are given by:
  \begin{equation}
    \Opt{\Pot}(\Angles)=-|\Angles|
    \ , \ 
    \OptTdens(\Angles)=\pi-|\Angles|
    \ , \ 
    \Angles\in \Sphere{1}
    .
  \end{equation}
  The graphs of $\OptTdens$ and the transported measures are
  reported in \cref{fig:circle}. The \OTD\ $\OptTdens$ progressively
  decays moving away from the point $\CutPoint$, until the antipodal
  point is reached, where it vanishes. According to the definition
  given in~\cite{Ambrosio:2003}, in this example the arcs $]-\pi,0[$
      and $]0,\pi[$ are the transport rays.
\end{example}

The extension of this idea to a general manifold leads \added{to} the following theorem,
which is a more complete version of~\cref{thm:opts-intro}.
\begin{theorem}
  \label{thm:opts}
  Let $(\Manifold,\Metric)$ be a compact and geodesically complete
  Riemannian manifold of dimension $\Mdim$ and with no boundary.
  Consider a point $\CutPoint\in\Manifold$ and the measures
  \begin{align}
    \label{eq:cut-source}
    \Source=\dv(\Manifold)\Dirac{\CutPoint}
    \  , \ 
    \Sink=\dv \; .
  \end{align}
  Let $(\OptPot,\OptTdens)$ be the solution of the
  \MKEQS\ with measures given by~\cref{eq:cut-source}. Then :
  \begin{enumerate}[i)]
  \item \label[stat]{enum:dist}
    the Kantorovich Potential $\OptPot$ coincides with
    minus the geodesic distance from $\CutPoint$:
    \begin{equation}
      \label{eq:optpot-cut}
      \OptPot(x)=-\Dist_{\Metric}(x,\CutPoint) \; ;
    \end{equation}
  \item \label[stat]{enum:continuity} when restricted to
    $\Manifold\setminus\left\{\CutPoint\right\}$, the
    \OTD\ $\OptTdens$ admits a \emph{continuous} density $\Density$
    with respect to the volume form $\dv$:
    \begin{equation*}
      \Opt{\Tdens} = \Density \dv \ , \
      \Density\in\Cont[0](\Manifold\setminus\left\{\CutPoint\right\}) \;;
    \end{equation*}
  \item  \label[stat]{enum:cutlocus}
    the zero set of $\Density$ coincides with the cut locus
    $\Cut{\CutPoint}$ of $\CutPoint$:
    \begin{equation*}
      \Cut{\CutPoint}
      =
      \left\{
        \Point\in\Manifold\setminus\left\{\CutPoint\right\}
        \ : \ \Density(\Point) = 0  
      \right\}\; ;
    \end{equation*}
  \item \label[stat]{enum:density} for any $\Metric[\CutPoint]$-orthonormal
    basis $\Base$, the local representative
    $\IntSetTdens^{\Base}:=\Density \circ (\Tilde{\Chart})^{-1}$ in
    the chart $\Chart^{\Base}$ of the function $\Density$ is given by
    \begin{equation}
      \label{eq:opttdens-cut}
      \begin{aligned}
        \IntSetTdens^{\Base}(\SphereDist,\SphereAngle)
        &=
        \frac{1}{ G^{\Base}(\SphereDist,\SphereAngle) \SphereDist^{\Mdim-1}} 
        \int_{\SphereDist}^{\TanCutTime_{\CutPoint}(\SphereAngle)}
        G^{\Base}(s,\SphereAngle) s^{\Mdim-1}\ds        
        \; ,
      \end{aligned}
    \end{equation}
    for all $(\SphereDist,\SphereAngle) \in
    \Tilde{\Chart}(\IntSet{\CutPoint}^{\Base})\subset \REALP\times
    \SphereAngles$.
  \end{enumerate}
\end{theorem}

\begin{proof}
  We first observe that thanks to~\cref{lemma:charts}, given two
  $\Metric[\CutPoint]$-orthonormal basis $\Base$ and $\Base'$ the
  functions $\IntSetTdens^{\Base}\circ \Chart^{\Base}$ and
  $\IntSetTdens^{\Base'}\circ \Chart^{\Base'}$ coincide in the
  intersection of their domains. Therefore, there exists a
  \emph{continuous} function $\Density$ on $\IntSet{\CutPoint}$ with
  local representative $\IntSetTdens^{\Base}$ in any chart
  $\Chart^{\Base}$ given by~\cref{eq:opttdens-cut}. Moreover, since
  for every $\Base$ $\IntSetTdens^{\Base}$ is positive in its domain
  $\Chart^{\Base}(\IntSet{\CutPoint}^{\Base})$ and, for every
  $\SphereAngle$,
  \begin{equation*}
    \lim_{\SphereDist \to \TanCutTime_{\CutPoint}(\SphereAngle)}
    \IntSetTdens^{\Base} (\SphereDist , \SphereAngle ) = 0 \; ,
  \end{equation*}
  the function $\Density$ admits a (unique) continuous extension on
  $\IntSet{\CutPoint}\cup\Cut{\CutPoint}=\Manifold\setminus\left\{\CutPoint\right\}$,
  which we will continue to denote with $\Density$, which is
  non-negative and whose zero set coincides with the cut locus of
  $\CutPoint$. This proves statements~\ref{enum:cutlocus}. The term
  $G^{\Base}$, defined in~\cref{eq:metric-in-coordinate}, is clearly
  bounded from above and from below. Hence, the integral
  in~\cref{eq:opttdens-cut} tends to zero as
  $\SphereDist \to \TanCutTime_{\CutPoint}(\SphereAngle)$.
    
  Now, we note that the source term $\Sink$ in~\cref{eq:cut-source} is
  absolutely continuous with respect to the volume form $\dv$.  Thus,
  \cref{lemma:prop-tdens} ensures that the \OTD\ $\OptTdens$ admits a
  unique density with respect to the volume form $\dv$. We now show
  that, with $\Density$ defined above, the pair $(\OptPot,\OptTdens)$
  given by
  \begin{equation}
    \label{eq:sol-pair}
    \OptPot(\Point) = -\Dist_{\Metric}(\Point,\CutPoint)
    \ ,\ \OptTdens\overset{a.e.}{=} \Density \dv
  \end{equation}
  solves~\cref{eq:mkeqs-manifold}.  It is clear that $\OptPot$
  satisfies the constraints
  in~\cref{eq:mkeqs-manifold-eikonal,eq:mkeqs-manifold-gradless1}.
  Thus, we only have to prove that $(\OptPot,\OptTdens)$ solves
  \cref{eq:mkeqs-weak-manifold} for $\Source $ and $\Sink$
  in~\cref{eq:cut-source} i.e.,
  \begin{equation}
    \label{eq:elliptic-cutlocus}
    \int_{\Manifold} 
    \Scal{\Grad[\Metric] \Opt{\Pot}}{\Grad[\Metric] \Test}_{ \Metric }
    d\OptTdens
    =    
    \int_{\Manifold}(\dv(\Manifold)) \Test \Dirac{\CutPoint}
    - 
    \int_{\Manifold}\Test \dv 
    \qquad 
    \forall\Test\in\Cont[1](\Manifold).
  \end{equation}
  For all $\Test\in\Cont[1](\Manifold)$, the right-hand side $RHS$
  of the above equation can be evaluated as
  \begin{equation}
    \label{eq:right-hand}
    RHS=(\dv(\Manifold))\Test(\CutPoint)
    -
    \int_{\Manifold}\!\Test\dv
    \; .
  \end{equation}
  We know, from \cref{lemma:cut-locus-prop} and~\cref{remark:line},
  that for any $\Metric[\CutPoint]$-orthonormal basis $\Base$ the set
  $\Manifold \setminus\IntSet{\CutPoint}^{\Base}$ has zero measure on
  $\Manifold$. Thus, we can restrict the integrals on the left-hand
  side of~\cref{eq:elliptic-cutlocus} to $\IntSet{\CutPoint}^{\Base}$
  and use the Riemannian polar coordinates to compute them. Write
  $\Tilde{\Test}=\Test\circ(\Tilde{\Chart})^{-1}$ and
  $\Tilde{\Pot}=\OptPot\circ(\Tilde{\Chart})^{-1}$. Thus,
  $\Tilde{\Pot}(\SphereDist,\SphereAngle)=-\SphereDist[]$ and
  the left-hand side
  $LHS$ of~\cref{eq:elliptic-cutlocus} becomes:
  \begin{multline*}
    LHS=\intsphere
    \int_{0}^{\TanCutTime_{\CutPoint} (\SphereAngle)}   
      \Grad[]
      \Tilde{\Pot}(\SphereDist[\ \!\!],\SphereAngle[\ \!\!])
    \cdot
    \Tilde{\Metric}(\SphereDist[\ \!\!],\SphereAngle[\ \!\!])^{-1}
    \Grad[]
      \Tilde{\Test}(\SphereDist[\ \!\!],\SphereAngle[\ \!\!]) 
      \Tilde{\IntSetTdens}(\SphereDist[],\SphereAngle[]) 
    \Tilde{\JacFull}(\SphereDist[\ \!\!],\SphereAngle[\ \!\!])
    d \SphereDist[\ \!\!] d \SphereAngle[\ \!\!] \; .
  \end{multline*}
  Now,
  $\Grad[]
  \Tilde{\Pot}(\SphereDist[\ \!\!],\SphereAngle[\ \!\!])
  =(-1,0,\ldots,0)$ and from~\cref{lemma:metric-polar-coordinates} we
  have
  \begin{equation*}
    \Grad[]
    \Tilde{\Pot}(\SphereDist[\ \!\!],\SphereAngle[\ \!\!])
    \cdot
    \Tilde{\Metric}(\SphereDist[\ \!\!],\SphereAngle[\ \!\!])
    ^{-1}
    \Grad[]
    \Tilde{\Test}(\SphereDist[\ \!\!],\SphereAngle[\ \!\!]) 
    =
    -\partial_{\SphereDist[\ \!\!]}
    \Tilde{\Test}(\SphereDist[\ \!\!],\SphereAngle[\ \!\!]) \; .
  \end{equation*}
  Thus, recalling~\cref{eq:opttdens-cut}
  and~\cref{eq:polar-volume-form} of~\cref{lemma:factorization}
  \begin{equation*}
    LHS=-
    \intsphere
    \int_{0}^{\TanCutTime_{\CutPoint} (\SphereAngle[\ \!\!])} 
    \left(
      \int_{\SphereDist[\ \!\!]}^{\TanCutTime_{\CutPoint} (\SphereAngle[\ \!\!])}
      \Tilde{\JacFull}(s,\SphereAngle[\ \!\!])
      d s
    \right)
    \partial_{\SphereDist[\ \!\!]}
    \Tilde{\Test}(\SphereDist[\ \!\!],\SphereAngle[\ \!\!]) 
    d \SphereDist[\ \!\!] d \SphereAngle[\ \!\!] \; .
  \end{equation*}
  Integration by parts yields:
  \begin{align*}
    \nonumber
    LHS&=
    -
    \intsphere
    \left(
      \left[
        \left(
          \int_{\SphereDist[\ \!\!]}^{\TanCutTime_{\CutPoint}(\SphereAngle[\ \!\!])}
          \Tilde{\JacFull}(s,\SphereAngle[\ \!\!])
          d s 
        \right)
        \Tilde{\Test}(\SphereDist[\ \!\!],\SphereAngle[\ \!\!]) 
      \right]_{0}^{\TanCutTime_{\CutPoint}(\SphereAngle[\ \!\!])}
    \right)d \SphereAngle[\ \!\!]
    \nonumber \\
    &\qquad
    +
    \intsphere
    \left(
      \int_{0}^{\TanCutTime_{\CutPoint}(\SphereAngle[\ \!\!])}
      \partial_{\SphereDist[\ \!\!]}\left(
        \int_{\SphereDist[\ \!\!]}^{\TanCutTime_{\CutPoint}(\SphereAngle[\ \!\!])}
        \Tilde{\JacFull}(s,\SphereAngle[\ \!\!])
        d s 
      \right)
      \Tilde{\Test}(\SphereDist[\ \!\!],\SphereAngle[\ \!\!]) 
      \
      d \SphereDist[\ \!\!]
    \right)
    d\SphereAngle[\ \!\!]
    \; .
  \end{align*}
  The first term is evaluated by taking the separate limits as
  $\SphereDist[\ \!\!]\Tendsto
  \TanCutTime_{\CutPoint}(\SphereAngle[\ \!\!])$ and as
  $\SphereDist[\ \!\!]\Tendsto 0$, with the former yielding zero. In conclusion,
  \begin{align*}
    LHS&=
      \Test(\CutPoint)
      \intsphere
      \int_{0}^{\TanCutTime_{\CutPoint}(\SphereAngle)} 
      \Tilde{\JacFull}(s,\SphereAngle)
      d s
      d \SphereAngle[\ \!\!]
         -
      \intsphere
      \int_{0}^{\TanCutTime_{\CutPoint}(\SphereAngle[\ \!\!])} 
      \Tilde{\JacFull}(\SphereDist[\ \!\!],\SphereAngle[\ \!\!])
      \Tilde{\Test}(\SphereDist[\ \!\!],\SphereAngle[\ \!\!]) 
      \
      d \SphereDist 
      d \SphereAngle
    \\
    &=(\dv(\Manifold))\Test(\CutPoint)
    -
    \int_{\IntSet{\CutPoint}^{\Base}}
    \Test
    \dv \; ,
  \end{align*}
  and then $LHS=RHS$ for all $\Test\in\Cont[1](\Manifold)$,
  proving~\cref{eq:elliptic-cutlocus}. This shows that the pair
  $(\OptPot,\OptTdens)$ in~\cref{eq:sol-pair} solves of the
  \MKEQS\ for $\Source $ and $\Sink$ in~\cref{eq:cut-source}.  This
  proves statements~\ref{enum:dist},\ref{enum:continuity}, and \ref{enum:density}.
  Statement~\ref{enum:cutlocus} has already been proved.
\end{proof}

\begin{Remark}
  It is worth noting that, thanks to~\cref{thm:opts} and to the
  results in~\cite{Facca-et-al-numeric:2020}, we can give a
  variational characterization of the \OTD\ as the minimizer of:
  \begin{equation*}
    \min_{\Tdens}
    \left\{\Lyap_{\Metric}({\Tdens}) \; : \; \Tdens \in
    \Cont[0](\Manifold\setminus\left\{\CutPoint\right\},\REALP) \right\}
  \end{equation*}
  with $\Lyap_{\Metric}=\Ene_{\Metric}+\Mass_{\Metric}$ and
  \begin{align*}
    \Ene_{\Metric}({\Tdens})
    &:= \sup_{\Test \in \Lip{}(\Manifold)}
    \left\{
    \int_{\Manifold}
    (\dv(\Manifold)\Dirac{\CutPoint} - \dv)\Test
    -
    {\Tdens}\frac{\ABS{\Grad[\Metric]\Test}^2}{2}
    \dv
    \right\}
    \nonumber \\
    \Wmass_{\Metric}({\Tdens})
    &:=\frac{1}{2}\int_{\Manifold}{{\Tdens}\dv} \;.
    \nonumber
  \end{align*}
  This provides a variational characterization of the cut locus of a
  point $\CutPoint\in\Manifold$ as the zero set of the minimizer
  $\OptTdens$ of $\Lyap_{\Metric}({\Tdens})$.
\end{Remark}

\section{Numerical approximation of the cut locus}
\label{sec:cut locus-experiments}
In this section we present our strategy for the calculation of the cut
locus based on~\cref{thm:opts}. Unfortunately, it is numerically
challenging to look directly at the zero set of the expression for
\OTD\ $\OptTdens$ given in~\cref{eq:opttdens-cut}. Indeed, this would
require the approximation of the distance and cut time functions,
i.e. the same unknowns in the identification problem of the cut locus.
As an alternative, at least in the case of a surface $\Surf$ embedded
in $\REAL^3$, we approach the numerical solution of the
\MKEQS~\cref{eq:mkeqs-manifold} by means of the ``dynamic''
reformulation of the \MKEQS, called \DMK, recently proposed
in~\cite{Facca-et-al:2018}, and its finite-element-based
discretization, described in~\cite{Facca-et-al-numeric:2020}. More
precisely, we use the extension to the surface setting of the \DMK\
approach as described in~\cite{Berti-et-al-numerics:2021}, in which
the numerical schemes developed in~\cite{Facca-et-al-numeric:2020} are
extended to the surface setting using the Surface Finite Element
Method (SFEM) framework reviewed in~\cite{Dziuk-Elliott:2013}.

We summarize here the fundamental steps of SFEM-DMK that impact on our
goal of calculating the cut locus.  First, the surface $\Surf$ is
decomposed with a geodesic triangulation $\Surf=\Triang(\Surf)$,
formed by triangles whose edges are the
geodesics between the vertices. Next, this triangulation is approximated
by its piecewise linear interpolant
$\Surfh=\TriangH(\Surf)=\cup \Cell[\Itdcell]$, i.e., the union
of 2-simplices $\Cell[\Itdcell]$ in $\REAL^3$ having the same
vertices as $\Triang(\Surf)$ (for more details
see~\cite{Morvan,Dziuk-Elliott:2013}). Using $\Surfh$ it is possible
to define appropriate discrete geometric quantities, such as surface
gradients and discrete finite element function spaces, that allow the
numerical discretization of the \MKEQS\ on surfaces embedded in
$\REAL^3$.  We refer to~\cite{Dziuk-Elliott:2013} for details on SFEM
and to~\cite{Berti-et-al-numerics:2021} for details and application
examples of SFEM to the \DMK\ on surfaces.

In practice, the calculation of the cut locus takes place on the
piecewise linear interpolation $\Surfh$ of $\Surf$ by means of the
SFEM discretization of the \MKEQS\ with
$\Source =\dv(\Surfh)\Dirac{\CutPoint}$ and $\Sink = \dv$.  To
maintain stability we follow the procedure described
in~\cite{Facca-et-al-numeric:2020} whereby the Kantorovich potential
$\OptPot$ is interpolated on a uniformly refined mesh
$\Triang[\MeshPar/2](\Surfh)\subset\TriangH(\Surfh)$ by piecewise
linear polynomials, while the OT density $\OptTdens$ is approximated
on $\TriangH(\Surfh)$ by piecewise constants.  As a consequence, the
approximated cut locus is formed by the union of the triangles
$\Cell[\Itdcell]\in\TriangH(\Surfh)$ where our numerically evaluated \OTD,
$\OptTdensH$, is close to zero.  Thus, we give the following
definition of the approximate cut locus ($\ACL$):
\begin{equation*}
  \ACL := 
  \left\{
    \cup \Cell[\Itdcell] \in \TriangH(\Surfh)
    \mbox{ s.t. } 
    (\OptTdensH )|_{\Cell[\Itdcell]} \le \TolCut \right\}
\end{equation*}
where $\TolCut$ is an appropriate preselected tolerance.
\begin{remark}
  The definition of this tolerance is crucial for a proper
  approximation of $\Cut{\CutPoint}$. Indeed, it is rather difficult
  to approximate a one-dimensional structure, and more so single
  points, as a union of triangles of $\TriangH(\Surf)$, and we will
  see by experimentation that different values of $\TolCut$ lead to
  different $\ACL$. Another difficulty arises from the use of the
  singular source term $\Source=\dv(\Surf)\Dirac{\CutPoint}$, which
  does not belong to the dual of $\Hilb{1}(\Surf)$. This is a typical
  problem in applications and yields suboptimal SFEM convergence rates
  for $\OptPot$. To address this issue one may either use
  regularization along the lines of~\cite{Tornberg:2003,Tornberg:2004}
  or a posteriori error estimations to adapt the surface mesh. We
  choose not to employ any of these approaches since first order
  global accuracy is in any case enforced by the piecewise constant
  discretization of $\OptTdens$.
\end{remark}

We would like to note that according to~\cite{ALBANO201651} the
cut locus is stable under $\mathcal{C}^2$ perturbations of the domain
or of the metric but stability may be lost in case of $\mathcal{C}^1$
perturbations. However, the ensuing experimental results show that
our piecewise constant approximation of the transport density
introduces enough regularization to our numerical solution, at least
for the sample problems addressed in this work. 

\subsection{Numerical experiments}
\label{sec:num-exp}

\added{First we would like to note that a Python notebook reproducing
  the experiments presented in this section can be found at this link
  \href{https://doi.org/10.5281/zenodo.5710660}{https://doi.org/10.5281/zenodo.5710660}.
  The source code is available at the
  following
  \href{https://gitlab.com/enrico\_facca/dmk_solver}{https://gitlab.com/enrico\_facca/dmk\_solver} (subdirectory FaccaBertiFassoPutti2021\_CutLocus).}
  %

\begin{figure}
  \centerline{
    \includegraphics[width=0.75\textwidth,
    ]{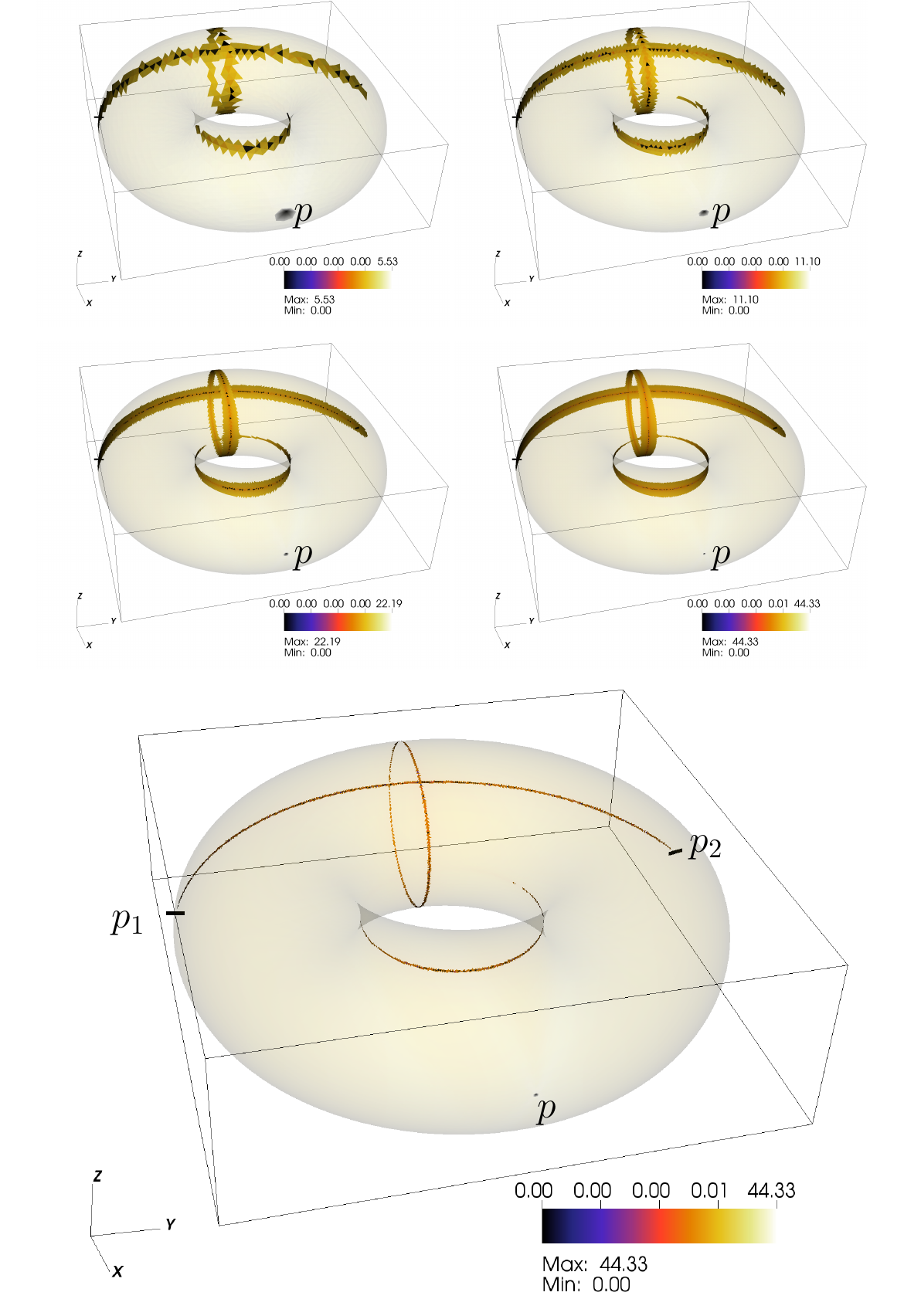}
  }
  \caption{ \added{Numerical identification of the $\ACL$ of a torus
      obtained at four uniform mesh refinement levels for
      $\TolCut=0.001$. The color map refers to $\OptTdensH$ values.
      The coarsest mesh level contains 2712 nodes and 5424 triangles.
      The two segments in black intersect the surface of the torus at
      points $p_1$ and $p_2$ described in \cref{fig:torus-chart}.  The
      bottom panel shows the $\ACL$ for $\TolCut=0.0001$ on the finest
      grid level.}}
  \label{fig:cut locus-torus}
\end{figure} 

The first test case concerns the example described in~\cref{sec:cut
  locus} and shown in~\cref{fig:torus-chart} and deals with the
calculation of the cut locus $\Cut{\CutPoint}$ of a torus $\Torus$
with radii $\Rmax=2$ and $\Rmin=1$ with respect to the point
$\CutPoint=(3,0,0)$.  The results of our simulations are shown
in~\cref{fig:cut locus-torus}.
The top panels of the figure depict the
$\ACL$ obtained using a tolerance $\TolCut=10^{-3}$ on four refinement
levels. It is clear that $\ACL$ approximates the real cut locus
with satisfactory accuracy already on the coarsest grid, with
progressively improving resolution as expected by the higher
refinement levels. The effect of the mesh finite size is clearly
discernible but no instabilities in the identification of the cut
locus are visible. Note that perturbations of $\Surf$ and $\Metric$ in
the refinement step from $\Triang[\MeshPar](\Surf)$ to
$\Triang[\MeshPar/2](\Surf)$ are not strictly
$\mathcal{C}^{1,1}$-regular. However, we can view $\TriangH(\Surf)$ as
a linear interpolation of $\Surf$ whose error can be bounded by
$\MeshPar^2$ times the norm of the second fundamental form of
$\Surf$~\cite{Morvan,Dziuk-Elliott:2013,art:BP20}. This regularity is
sufficient to provide empirical justification of the stability of our
calculations. This rationale is further strengthen\added{ed} by the results at
the finest grid level.

The $\ACL$ for $\TolCut=10^{-4}$ is shown in the bottom
panel of~\cref{fig:cut locus-torus}. A much better approximation of
the real cut locus is displayed with increased level of
details. However, one portion of the ``internal equator'' opposite to
$\CutPoint$ is not identified.
We can give a heuristic explanation for this phenomenon by looking at
the explicit formula of the \OTD\ in \cref{eq:opttdens-cut}.  Roughly,
this formula predicts an increased value of the \OTD\ in
the regions of larger mass fluxes.
In our case this corresponds to the region surrounding the internal
equator, where the geodesics starting at $\CutPoint$ arrive with a
small angle with respect to the internal great circle. Thus, the
values of $\OptTdensH$ in these triangles are relatively large, and
the approximation of the cut locus becomes problematic when using an
absolute identification criterion.  This problem could of course be relieved by
employing standard strategies that combine relative error measures
and adaptive mesh refinements, tasks that go beyond the purpose of this
paper.

The second test case is taken from~\cite[Section
  5]{Itoh-Sinclair:2004} where the authors consider the triaxial
ellipsoid given by
\begin{equation*}
  \Ellipsoid:=
  \left\{
    (x,y,z) \in \REAL^3 \ : \
    (x/0.2)^2+(y/0.6)^2+z^2=1 
  \right\}.
\end{equation*}
The cut loci of two points in $\Ellipsoid$ are studied. The first
point considered is $\CutPoint^{1}=(-0.115470, 0, 0.816497)$, an
umbilic point whose cut locus consists of a single point. The second
cut locus is relative to the point
$\CutPoint^{2}=(-0.151128, -0.350718,0.295520)$ and is conjectured to
be an arc on the opposite site of the ellipsoid~\cite[Conjecture
5.2]{Itoh-Sinclair:2004}.
\begin{figure}
  \centerline{
    \includegraphics[
      trim={0.9cm 1.5cm 6.5cm 2.3cm},clip,
      width=0.32\textwidth,
    ]{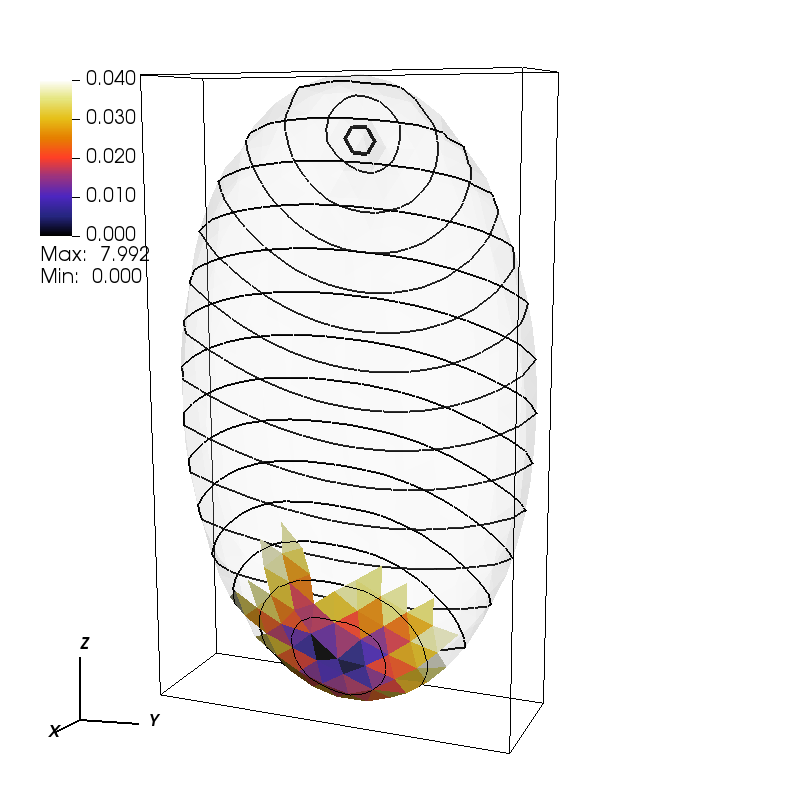}
    \includegraphics[
      trim={0.9cm 1.5cm 6.5cm 2.3cm},clip,
      width=0.32\textwidth,
    ]{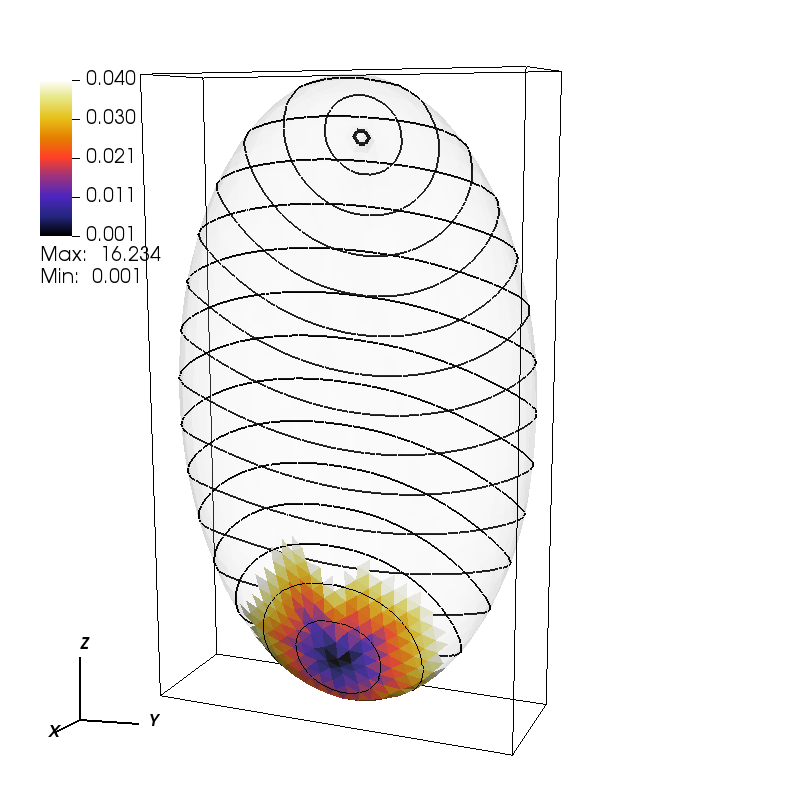}
    \includegraphics[
      trim={0.9cm 1.5cm 6.5cm 2.3cm},clip,
      width=0.32\textwidth,
    ]{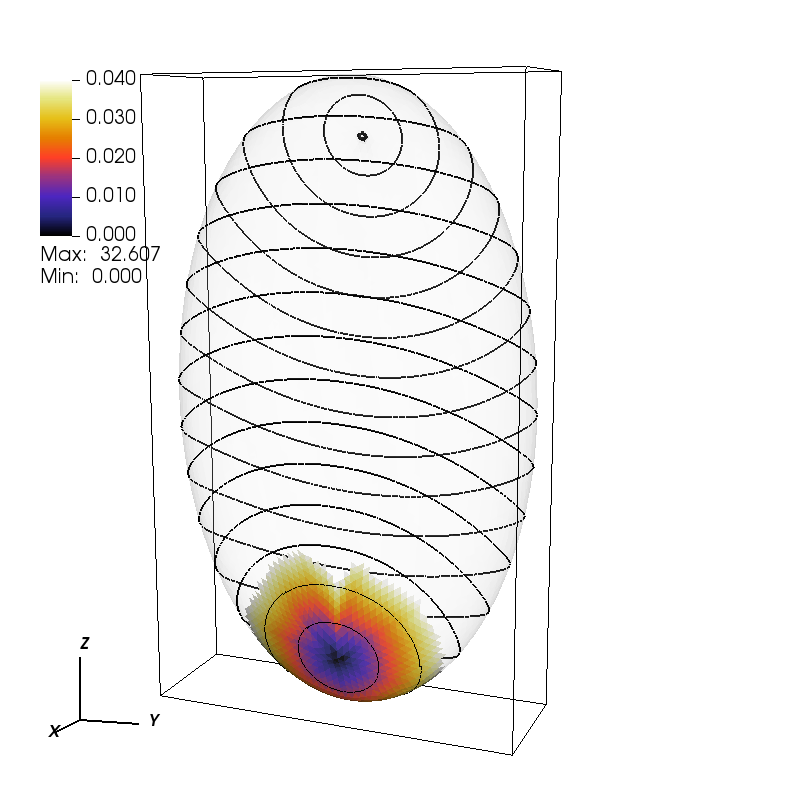}
  }
  \centerline{
    \includegraphics[
      trim={0.6cm 1.5cm 6.5cm 2cm},clip,
      width=0.32\textwidth,
    ]{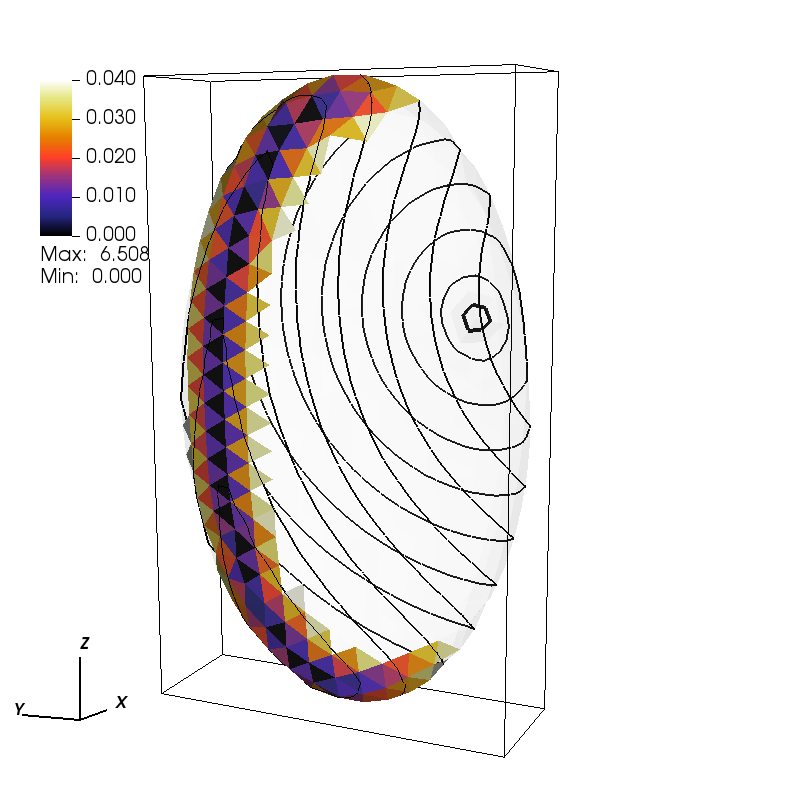}
    \includegraphics[
      trim={0.6cm 1.5cm 6.5cm 2.cm},clip,
      width=0.32\textwidth,
    ]{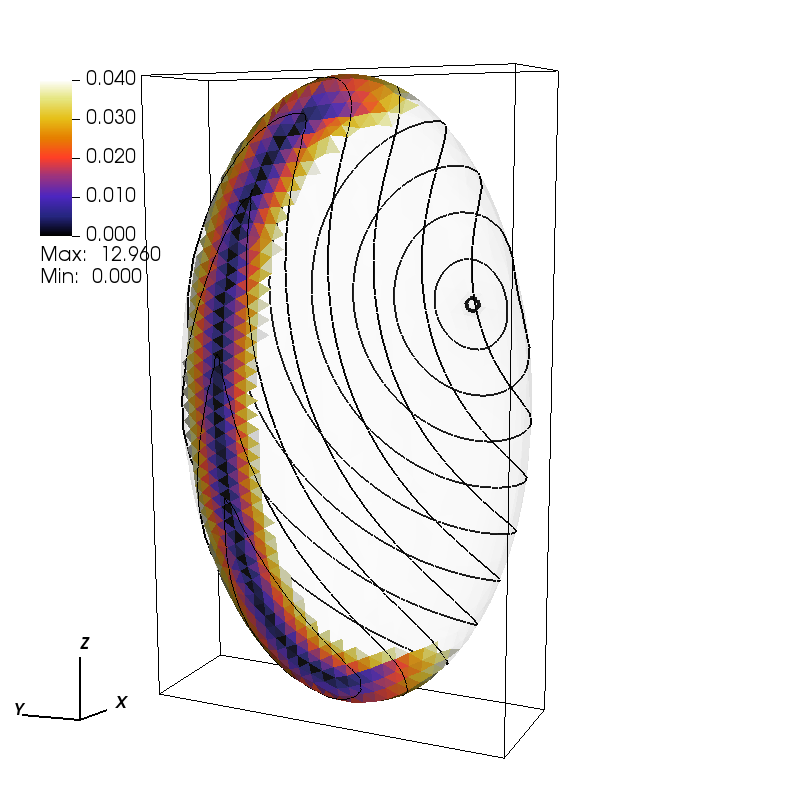}
    \includegraphics[
      trim={0.6cm 1.5cm 6.5cm 2.cm},clip,
      width=0.32\textwidth,
    ]{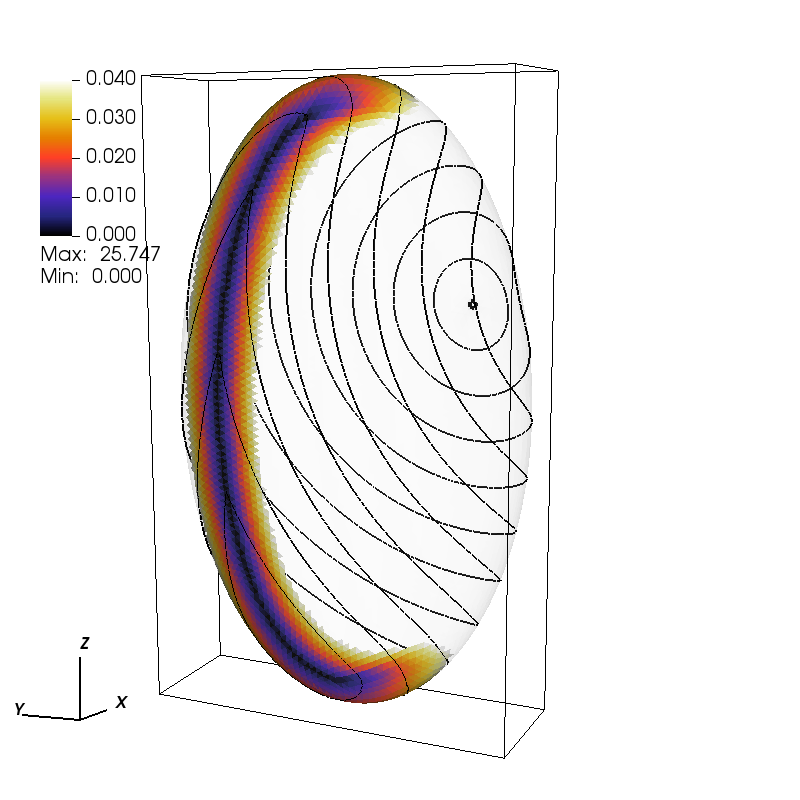}
  }
  \caption{ Ellipsoid $(x/0.2)^2+(y/0.6)^2+(z)^2=1$. The upper panels
    reports the spatial distribution of $\OptTdensH$ for the point
    $\CutPoint^{1}=(-0.115470, 0, 0.816497)$. We report the results for
    three meshes, one the conformal refinement of the other, with the
    first having 534 nodes and 1064 triangles, and the latter 8514
    nodes 17024 triangles. The point $\CutPoint^{1}$ is located
    ``behind'' the visible ellipsoid an\added{d} it is marked with a black
    circle. We report only those triangles where $\OptTdensH$ is below
    $0.04$. We also report the contour lines of $\OptPotH$. On the
    lower panels we report the same plots for the
    $\CutPoint^{2}=(-0.151128,-0.350718,0.295520)$.}
  \label{fig:cut locus-ellipsoid}
\end{figure} 
In~\cref{fig:cut locus-ellipsoid} we report the spatial distribution
of the approximate \OTD\ $\OptTdensH$ associated to the points
$\CutPoint^{1}$ and $\CutPoint^{2}$, on the top and bottom panels,
respectively.  We start from an initial triangulation
$\Triang[\MeshPar]$ of the ellipsoid with 534 nodes and 1064
triangles, obtained by means of the software described
in~\cite{Persson-Strang:2004}. We generate a sequence of finer grids
conformally refining $\Triang[\MeshPar]$ and ``lifting'' the added
nodes (the mid-points of the edges of the triangles) to the
ellipsoid. We report the spatial distribution of $\OptTdensH$ only on
those triangles where it is below the threshold $\OptTdensH<0.04$ to
appreciate the decay of the \OTD\ as it approaches the cut
locus. Already at the coarsest grid, the region were $\OptTdensH$
attains the lowest values (in black in~\cref{fig:cut locus-ellipsoid})
strongly resembles the approximate cut locus computed
in~\cite{Itoh-Sinclair:2004}.

\begin{figure}
  \centerline{
    \includegraphics[trim={0.9cm 1.5cm 6.5cm 0.cm},clip,
    width=0.49\textwidth]{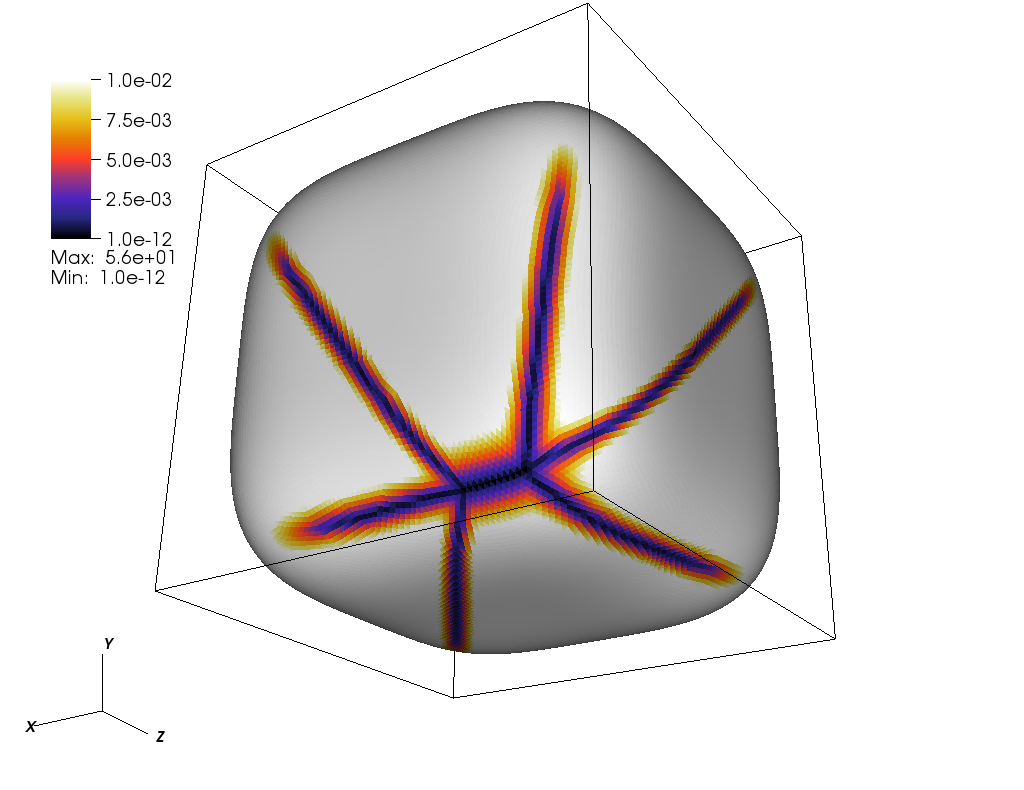}
    \quad
     \includegraphics[trim={0.9cm 1.5cm 6.5cm 0.cm},clip,
    width=0.49\textwidth]{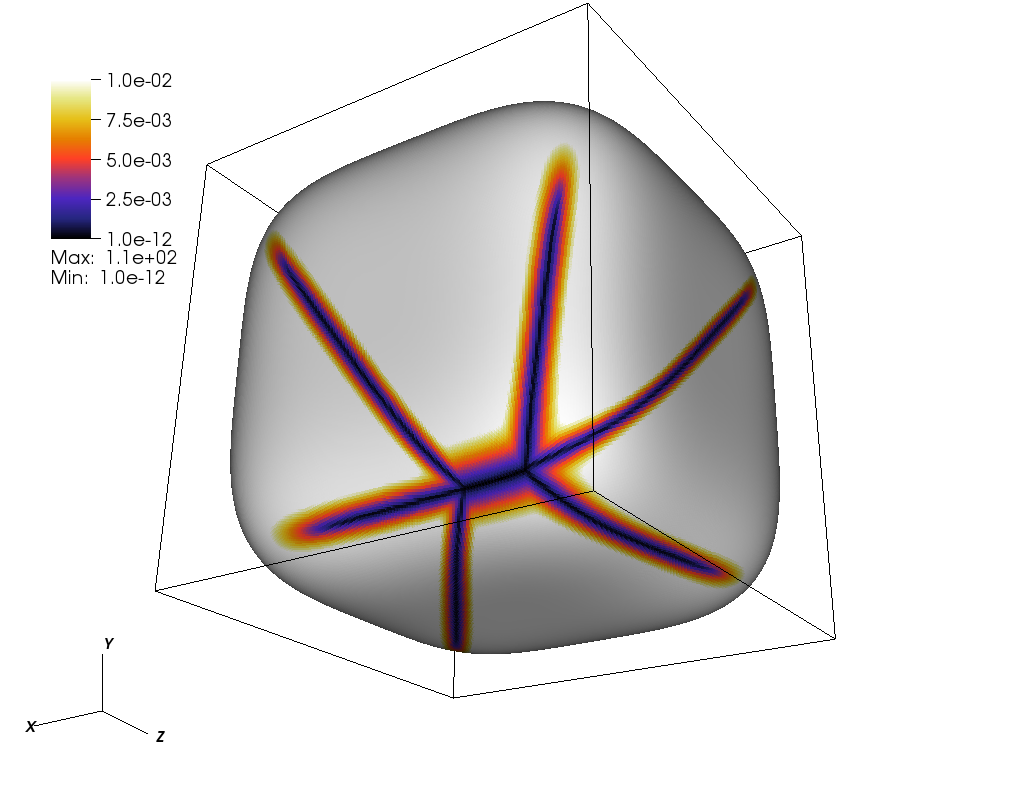}
  }
  \caption{ Spatial distribution of the $\OptTdensH$ on the
    surface $\Quartic$ using
    $\Source=|\Quartic|\Dirac{\CutPoint}$ with
    $\CutPoint=(0.533843,0.800764,0.844080)$ and
    $\Forcing=1$.  The left panel report the result using a
    mesh containing $34178$ nodes and $68352$ cells. The color
    scale has its maximum at $0.01$ to highlight the region
    where $\OptTdensH$ attains the lowest value.  The right
    panel reports the same results using the conformal
    refinement of the surface triangulation ($136706$ nodes
    and $273408$ cells). 
}
  \label{fig:cut locus-quartic}
\end{figure} 

The final test-case, again taken from~\cite{Itoh-Sinclair:2004}, looks
for the cut locus of the quartic equation defined by
\begin{equation*}
  \Quartic:=
  \left\{
    (x,y,z) \in \REAL^3 \ : \
    x^4+y^4+z^4=1
  \right\}
\end{equation*}
with respect to the point $\CutPoint=(0.533843,0.800764,0.844080)$.
\Cref{fig:cut locus-quartic} shows the spatial distribution of
$\OptTdensH$ obtained on two different triangulations. The first mesh
(left panel) is characterized by 34178 nodes and 68352 triangles in
$\TriangH$, where $\OptTdensH$ is defined, and 136712 nodes and 273408
triangles, in $\Trianghh$, where $\OptPotH$ lives.
The second mesh level (right panel) is exactly four times larger.  The
approximate zero set of $\OptTdensH$ shown on the two refinements
in~\cref{fig:cut locus-quartic} compares well with the approximate cut
locus reported in~\cite{Itoh-Sinclair:2004}.

We would like to note here that our DMK-based approach is much more
computationally efficient than ~\cite{Itoh-Sinclair:2004}. Indeed, the
computational cost for these simulations on a laptop computer equipped
with a 2014 Intel Core-I5 processor with 8Gbyte of RAM are 782 and
4606 seconds for the two meshes. This is to be compared with the
computational cost of 35506 seconds reported
in~\cite{Itoh-Sinclair:2004} to solve the same problem on a mesh with
49152 triangles (with an unspecified CPU). The difference in
performance has to be attributed to the fact that the algorithm
in~\cite{Itoh-Sinclair:2004} has a computational complexity that grows
exponentially with the size of the triangulation, while our algorithm
is affected by the classical polynomial computational complexity of
FEM methods.
Finally, we would like to note that the computational cost of our
method is comparable if not better to the approach described
in~\cite{generau2020numerical}. The use of implicit time-stepping in
combination with Newton method as proposed in~\cite{FaBe2020}, which
allows a drastic improvement in computational efficiency, is the next
step in our future studies.

\section{Conclusions}
We presented a new result showing the one-to-one correspondence
between the cut locus of a point $\CutPoint$ in a Riemannian manifold
$(\Manifold,\Metric)$ and the zero set of the \OTD\ $\OptTdens$
solution of the \MKEQS\ with $\Source=\dv(\Manifold)\Dirac{\CutPoint}$
and $\Sink=\dv$. \added{This new PDE-based characterization allows us
  to exploit standard finite element methods in combination with
  optimization techniques for the numerical approximation of the cut
  locus.}  Based on this result, we proposed a novel numerical
approach for the identification of the cut locus of a point on 2d
surfaces embedded in $\REAL^3$ using the \DMK\ method proposed
in~\cite{Berti-et-al-numerics:2021} for the solution of the of the
\MKEQS.  Numerical tests on few examples show that the cut locus can
be efficiently identified with the developed strategy.
The proposed \DMK-based numerical approach can be easily extended to
manifolds with dimension greater than two as long as 
the numerical solution of the PDE~\cref{eq:mkeqs-manifold-elliptic}
is \replaced{feasible}{ readily available}.

\section*{Acknowledgments}
\added{EF has been supported by the Centro di Ricerca Matematica ``Ennio De
Giorgi''.} FF has been partially supported by MIUR-PRIN project 20178CJA2B
\emph{New frontiers of Celestial Mechanics: theory and applications}.
MP has been partially supported by UniPD-SID-2016 project
\emph{Approximation and discretization of PDEs on Manifolds for
Environmental Modeling}.

\bibliographystyle{unsrt}  
\bibliography{strings,biblio_abbr,pubblication_ef}

\begin{thebibliography}{10}

\bibitem{chavel_2006}
Isaac Chavel.
\newblock {\em Riemannian Geometry: A Modern Introduction}.
\newblock Cambridge Studies in Advanced Mathematics. Cambridge University
  Press, 2 edition, 2006.

\bibitem{Mantegazza:2002}
Carlo Mantegazza and Andrea~Carlo Mennucci.
\newblock {H}amilton-{J}acobi equations and distance functions on {R}iemannian
  manifolds.
\newblock {\em Appl. Math. Optim.}, 47(1):1--25, January 2003.

\bibitem{Bonnard2009}
Bernard Bonnard, Jean~Baptiste Caillau, Robert Sinclair, and Minoru Tanaka.
\newblock Conjugate and cut loci of a two-sphere of revolution with application
  to optimal control.
\newblock {\em Ann. Inst. H. Poincar\'e Anal. Non Lin\'eaire},
  26(4):1081--1098, 2009.

\bibitem{Itoh-Sinclair:2004}
Jin~Ichi Itoh and Robert Sinclair.
\newblock Thaw: A tool for approximating cut loci on a triangulation of a
  surface.
\newblock {\em Experiment. Math.}, 13(3):309--325, 2004.

\bibitem{Misztal:2011}
Marek~Krzysztof Misztal, Jakob~Andreas B{\ae}rentzen, Francois Anton, and Steen
  Markvorsen.
\newblock Cut locus construction using deformable simplicial complexes.
\newblock In {\em 2011 Eighth International Symposium on Voronoi Diagrams in
  Science and Engineering (ISVD)}, pages 134--141. IEEE, June 2011.

\bibitem{Dey-Kuiyu:2009}
Tamal~K Dey and Kuiyu Li.
\newblock Cut locus and topology from surface point data.
\newblock In {\em Proceedings of the twenty-fifth annual symposium on
  Computational geometry}, pages 125--134, 2009.

\bibitem{Crane:2013}
Keenan Crane, C~Weischedel, and M~Wardetzky.
\newblock Geodesics in heat: {A} new approach to computing distance based on
  heat flow.
\newblock {\em ACM Trans. Graph.}, 2013.

\bibitem{generau2020cut}
Fran{\c{c}}ois G{\'e}n{\'e}rau, Edouard Oudet, and Bozhidar Velichkov.
\newblock Cut locus on compact manifolds and uniform semiconcavity estimates
  for a variational inequality.
\newblock arXiv preprint arXiv:2006.07222, 2020.

\bibitem{generau2020numerical}
Fran{\c{c}}ois G{\'e}n{\'e}rau, Edouard Oudet, and Bozhidar Velichkov.
\newblock Numerical computation of the cut locus via a variational
  approximation of the distance function.
\newblock arXiv preprint arXiv:2006.08240, 2020.

\bibitem{Figalli2008}
Alessio Figalli and C{\'e}dric Villani.
\newblock An approximation lemma about the cut locus, with applications in
  optimal transport theory.
\newblock {\em Methods and Applications of Analysis}, 15(2):149--154, June
  2008.

\bibitem{Figalli:2010}
Alessio Figalli and Ludovic Rifford.
\newblock Mass transportation on sub-riemannian manifolds.
\newblock {\em Geom. Funct, Anal.}, 20(1):124--159, Jun 2010.

\bibitem{Figalli-et-al:2011}
Alessio Figalli, Ludovic Rifford, and C{\'e}dric Villani.
\newblock Tangent cut loci on surfaces.
\newblock {\em Differ. Geom. Appl.}, 29(2):154--159, 2011.

\bibitem{Villani2011}
C{\'e}dric Villani.
\newblock Regularity of optimal transport and cut locus: From nonsmooth
  analysis to geometry to smooth analysis.
\newblock {\em Discrete Contin. Dyn. Syst. Ser. A}, 30(2):559--571, January
  2011.

\bibitem{Ambrosio:2003}
Luigi Ambrosio.
\newblock Lecture notes on optimal transport problems.
\newblock In {\em Lecture Notes in Mathematics}, pages 1--52. Springer, Berlin,
  Heidelberg, 2003.

\bibitem{Villani:2003}
C{\'e}dric Villani.
\newblock {\em Topics in optimal transportation}.
\newblock Number~58 in Graduate studies in mathematics. American Mathematical
  Soc., 2003.

\bibitem{Villani:2008}
C{\'e}dric Villani.
\newblock {\em Optimal Transport: Old and New}, volume 338.
\newblock Springer Science \& Business Media, 2008.

\bibitem{Santambrogio:2015}
Filippo Santambrogio.
\newblock {\em Optimal Transport for Applied Mathematicians}.
\newblock Birk{\"a}user, NY, 2015.

\bibitem{Pratelli:2005}
Aldo Pratelli.
\newblock Equivalence between some definitions for the optimal mass transport
  problem and for the transport density on manifolds.
\newblock {\em Ann. di Mat. Pura ed Appl.}, 184:215--238, 06 2005.

\bibitem{Evans-Gangbo:1999}
Lawrence~C Evans and Wilfrid Gangbo.
\newblock {\em Differential equations methods for the {M}onge--{K}antorovich
  mass transfer problem}, volume 653.
\newblock American Mathematical Soc., 1999.

\bibitem{Berti-et-al-numerics:2021}
Luca Berti, Enrico Facca, and Mario Putti.
\newblock Numerical solution of the {L}$^1$-optimal transport problem on
  surfaces.
\newblock ArXiv preprint available, 2021.

\bibitem{Facca-et-al-numeric:2020}
Enrico Facca, Sara Daneri, Franco Cardin, and Mario Putti.
\newblock Numerical solution of monge--kantorovich equations via a dynamic
  formulation.
\newblock {\em J. Scient. Comput.}, 82(3):1--26, 2020.

\bibitem{Facca-et-al:2018}
Enrico Facca, Franco Cardin, and Mario Putti.
\newblock Towards a stationary {M}onge--{K}antorovich dynamics: The physarum
  polycephalum experience.
\newblock {\em SIAM J. Appl. Math.}, 78(2):651--676, 2018.

\bibitem{Sakai1996}
Takashi Sakai.
\newblock {\em Riemannian geometry}, volume 149.
\newblock American Mathematical Soc., 1996.

\bibitem{Itoh-Minoru:2001}
Jin~Ichi Itoh and Minoru Tanaka.
\newblock The {L}ipschitz continuity of the distance function to the cut locus.
\newblock {\em Trans. Amer. Math. Soc.}, 353(1):21--40, 2001.

\bibitem{kimball1930}
Bradford~Fisher Kimball.
\newblock Geodesics on a toroid.
\newblock {\em Am. J. Math.}, 52(1):29--52, 1930.

\bibitem{Gravensen2005}
Jens Gravesen, Steen Markvorsen, Robert Sinclair, and Minoru Tanaka.
\newblock The cut locus of a torus of revolution.
\newblock {\em Asian J. Math.}, 9(1):103--120, 03 2005.

\bibitem{Jost:2008}
J{\"u}rgen Jost.
\newblock {\em Riemannian geometry and geometric analysis}, volume 42005.
\newblock Springer, 2008.

\bibitem{Blumenson:1960}
L.~E. Blumenson.
\newblock A derivation of n-dimensional spherical coordinates.
\newblock {\em Am. Math. Mon.}, 67(1):63--66, 1960.

\bibitem{Bouchitte-Buttazzo:2001}
Guy Bouchitt{\'e} and Giuseppe Buttazzo.
\newblock Characterization of optimal shapes and masses through
  {M}onge-{K}antorovich equation.
\newblock {\em J. Eur. Math. Soc.}, 3(2):139--168, 2001.

\bibitem{Feldman-McCann:2001}
Mikhail Feldman and Robert~J. McCann.
\newblock Monge's transport problem on a riemannian manifold.
\newblock {\em Trans. Amer. Math. Soc.}, 354(4):1667--1697, 2001.

\bibitem{Feldman-McCann:2002}
Mikhail Feldman and Robert~J. McCann.
\newblock Uniqueness and transport density in monge's mass transportation
  problem.
\newblock {\em Calc. Var. Partial Diff. Equ.}, 15(1):81--113, 2002.

\bibitem{DePascale-et-al:2004}
Luigi De~Pascale, Lawrence~C Evans, and Aldo Pratelli.
\newblock Integral estimates for transport densities.
\newblock {\em Bull. London Mat. Soc.}, 36(3):383--395, 2004.

\bibitem{Santambrogio:2009}
Filippo Santambrogio.
\newblock Absolute continuity and summability of transport densities: simpler
  proofs and new estimates.
\newblock {\em Calc. Var. Partial Diff. Equ.}, 36(3):343--354, 2009.

\bibitem{Dziuk-Elliott:2013}
Gerhard Dziuk and Charles~M Elliott.
\newblock Finite element methods for surface pdes.
\newblock {\em Acta Num.}, 22:289--396, 2013.

\bibitem{Morvan}
Jean-Marie Morvan.
\newblock {\em Generalized Curvatures}, volume~2 of {\em Geometry and
  Computing}.
\newblock Springer Science {\&} Business Media, Berlin, Heidelberg, May 2008.

\bibitem{Tornberg:2003}
Anna-Karin Tornberg and Bjorn Engquist.
\newblock Regularization techniques for numerical approximation of {PDE}s with
  singularities.
\newblock {\em J. Scient. Comput.}, 19(1-3):527--552, December 2003.

\bibitem{Tornberg:2004}
Anna-Karin Tornberg and Bjorn Engquist.
\newblock Numerical approximations of singular source terms in differential
  equations.
\newblock {\em J. Comp. Phys.}, 200(2):462--488, November 2004.

\bibitem{ALBANO201651}
Paolo Albano.
\newblock On the stability of the cut locus.
\newblock {\em Nonlinear Anal. Theory Methods Appl.}, 136:51 -- 61, 2016.

\bibitem{art:BP20}
Elena Bachini and Mario Putti.
\newblock Geometrically intrinsic modeling of shallow water flows.
\newblock {\em Math. Model. Num. Anal.}, 54(6):2125--2157, 2020.

\bibitem{Persson-Strang:2004}
Per-Olof Persson and Gilbert Strang.
\newblock A simple mesh generator in {M}atlab.
\newblock {\em SIAM Rev.}, 46(2):329--345, 2004.

\bibitem{FaBe2020}
Enrico Facca and Michele Benzi.
\newblock Fast iterative solution of the optimal transport problem on graphs.
\newblock {\em SIAM J. Sci. Comput.}, 2021.
\newblock To appear, ArXiv preprint available.

\end{thebibliography}

\end{document}